\documentclass[preprint,11pt,numbers,nopreprintline]{elsarticle}

\usepackage[T1]{fontenc}
\usepackage[utf8]{inputenc}
\usepackage{lmodern}
\usepackage{amsmath,amssymb,amsthm,mathtools,mathrsfs}
\usepackage{microtype,needspace}
\usepackage[a4paper,top=26mm,bottom=26mm,left=27mm,right=27mm]{geometry}
\usepackage[hidelinks]{hyperref}
\biboptions{sort&compress}

\allowdisplaybreaks[3]
\numberwithin{equation}{section}
\makeatletter
\let\ps@pprintTitle\ps@plain
\makeatother

\newtheorem{theorem}{Theorem}[section]
\newtheorem{proposition}[theorem]{Proposition}
\newtheorem{lemma}[theorem]{Lemma}

\theoremstyle{definition}

\theoremstyle{remark}
\newtheorem{remark}[theorem]{Remark}

\DeclareMathOperator{\tr}{tr}
\DeclareMathOperator{\trace}{trace}

\hypersetup{
 pdftitle={Biharmonic submanifolds in nonpositive space forms},
 pdfauthor={Yu Fu, Min-Chun Hong and Dan Yang},
 pdfsubject={Biharmonic submanifolds, PNMCV, and flat normal bundle},
 pdfkeywords={Biharmonic maps, Generalized Chen's conjecture, PNMCV, Flat normal bundle}
}
\journal{Journal of Functional Analysis}

\begin{document}

\begin{frontmatter}

\title{Biharmonic submanifolds in nonpositive space forms}
\author[dufe]{Yu Fu}
\ead{yufu@dufe.edu.cn}
\author[uq]{Min-Chun Hong}
\ead{hong@maths.uq.edu.au}
\author[lnu]{Dan Yang}
\ead{dlutyangdan@126.com}
\address[dufe]{Center of Applied Mathematics, Dongbei University of Finance and Economics,\\
Dalian 116025, P.\ R.\ China}
\address[uq]{Department of Mathematics, The University of Queensland,\\
Brisbane, QLD 4072, Australia}
\address[lnu]{School of Mathematics and Statistics, Liaoning University,\\
Shenyang 110036, P.\ R.\ China}

\begin{abstract}
We prove that there are no proper biharmonic submanifolds with parallel normalized
mean curvature vector field  in a nonpositive space form,
in arbitrary dimension and codimension. We also prove that every
biharmonic surface with flat normal bundle in a four-dimensional
nonpositive space form is minimal. These results provide substantial
support for Chen’s conjecture and the generalized Chen’s conjecture for hyperbolic space
in higher codimensions.
\end{abstract}


\end{frontmatter}

\section{Introduction}

Biharmonic maps are critical points of the bienergy functional, which
is given by
\[
E_2(\varphi,\Omega)
=\frac12\int_\Omega |\tau(\varphi)|^2\,\mathrm d v_g,
\qquad \tau(\varphi)=\trace_g\nabla\mathrm d\varphi,
\]
with respect to variations supported in relatively compact domains
$\Omega\subset M$. Jiang~\cite{Jiang1986} derived the corresponding
Euler--Lagrange equation, which is given by $\tau_2(\varphi)=0$, where
$\tau_2(\varphi)$ denotes the bitension field. For an isometric immersion
$\varphi\colon M^m\to N^{m+p}$, the tension field is given by
$\tau(\varphi)=mH$, where $H$ denotes the mean curvature vector field.
Since every minimal immersion is biharmonic, the classification problem
concerns nonminimal biharmonic immersions, which are called
\emph{proper biharmonic}. Various fundamental analytic properties of biharmonic maps have been investigated in the literature. Unique continuation theorems for biharmonic maps were established in~\cite{BrandingOniciuc2019}.

For immersions into Euclidean space, biharmonicity means that each
Cartesian component of $H$ is harmonic. The following longstanding conjecture was 
formulated by B.-Y. Chen in 1991~\cite{Chen1991}.

\medskip
\noindent\textbf{Chen's conjecture.}
\textit{Every biharmonic submanifold in a Euclidean space is minimal.}
\medskip

In 1992, Dimitri\'c~\cite{Dimitric1992} proved minimality for biharmonic
hypersurfaces with at most two distinct principal curvatures and for
pseudo-umbilical biharmonic submanifolds of dimension different from
four. For four-dimensional pseudo-umbilical submanifolds, the tangential
biharmonic equation is identically satisfied. In codimension two,
Maeta~\cite{MaetaPseudo2026} recently proved that every pseudo-umbilical
biharmonic submanifold $M^4$ in $N^6(c)$ with $c\leq0$ is minimal.

The hypersurface case in $\mathbb R^4$ was settled 
in~\cite{HV1995}, and subsequently in $\mathbb R^5$ and
$\mathbb R^6$ in \cite{FHZ2021,FHZ2023},
respectively.  More recently, the hypersurface case was established
in all dimensions~\cite{Abraham2026,FHY2026}.
In higher codimension, minimality was proved under the global
assumption that the immersion is proper~\cite{AM2013}.
For the related evolution problem, Fu, Hong and
Tian~\cite{FHT2026} obtained estimates for the maximal existence
time of the biharmonic
hypersurface flow in arbitrary dimension by combining
Gagliardo--Nirenberg inequalities with local energy estimates.

The generalized Chen's conjecture was proposed
in~\cite{CMO2001,CMO2002} for ambient manifolds of nonpositive
sectional curvature. Counterexamples have been constructed in
ambient manifolds with nonconstant negative sectional
curvature~\cite{OuTang2012}. For hyperbolic space, the conjecture
is stated as follows.

\medskip
\noindent\textbf{Generalized Chen's conjecture for hyperbolic space.}
\textit{Every biharmonic submanifold in a hyperbolic space is minimal.}
\medskip

The conjecture was proved for surfaces in $\mathbb H^3$
~\cite{CMO2002} and for hypersurfaces in $\mathbb H^4$
~\cite{BMO2010}. In arbitrary dimension, minimality was shown for
biharmonic hypersurfaces in hyperbolic space with at most two distinct
principal curvatures~\cite{BMO2008}. The hypersurface case was recently
settled in all dimensions by the authors~\cite{FHY2026}.

Progress in arbitrary codimension has also been made under global
assumptions. Nakauchi and Urakawa \cite{NU2013} proved that a complete
biharmonic submanifold of a Riemannian manifold with nonpositive
sectional curvature is minimal if the mean curvature vector field is
square integrable. Maeta~\cite{Maeta2014} proved minimality for properly
immersed biharmonic submanifolds of hyperbolic space.
Very recently, Maeta~\cite{MaetaWintgen2026} also proved that biharmonic Wintgen ideal
submanifolds in nonpositive space forms are minimal, provided that the
submanifold has dimension at least three or is a surface with codimension
at most two.

Thus Chen's conjecture and its hyperbolic-space analogue have been
settled for hypersurfaces. The remaining general problem concerns mainly submanifolds of higher codimension.
In this paper we consider submanifolds of higher codimension under geometric
conditions on the normal bundle.

Proper biharmonic submanifolds occur naturally in
spheres~\cite{CMO2001,CMO2002}. The conjecture of Balmu\c{s}, Montaldo
and Oniciuc~\cite{BMO2008} asserts that \emph{every biharmonic
	submanifold of a Euclidean sphere has constant mean curvature}.
The constant mean curvature case for proper biharmonic surfaces in
$\mathbb S^4$ was classified by Balmu\c{s} and
Oniciuc~\cite{BO2009}. Mean curvature gap and partial classification
results were obtained under the assumption of a parallel mean curvature
vector field~\cite{BO2012}, and PNMCV submanifolds were studied
in~\cite{BMO2013}. Bochner--Simons formulas also yielded rigidity
theorems for biharmonic submanifolds in space forms~\cite{FLO2021}.
For hypersurfaces in Euclidean spheres, unique continuation of constant mean
curvature was established under an additional gradient
condition~\cite{BLO2021}. Further results can be found
in~\cite{OuChen2020,FetcuOniciuc2022}.

Throughout the paper, $N^{m+p}(c)$ denotes a Riemannian space form of
constant sectional curvature $c$. We write $h$ for the second
fundamental form, $\nabla^\perp$ for the normal connection, and $A_\nu$
for the shape operator in the normal direction $\nu$.

\Needspace{9\baselineskip}
An immersion has a parallel normalized mean curvature vector field
(PNMCV) if $H$ is nowhere zero and the unit normal field
$\xi=H/|H|$ is parallel with respect to the normal connection.
Writing $f=|H|$, the PNMCV assumption can be expressed as
\[
f>0,\qquad \xi=\frac{H}{f},\qquad
\nabla^\perp\xi=0\quad\text{on }M.
\]
The function $f$ need not be constant. In fact, a PNMCV immersion has parallel
mean curvature vector field if and only if $f$ is locally constant.

Chen~\cite{Chen2019} proved the nonexistence of biharmonic PNMCV
surfaces in Euclidean spaces of arbitrary dimension and conjectured
the same conclusion in higher dimensions. The three-dimensional case
in $\mathbb R^5$ was subsequently established in~\cite{ST2023}.
Our first theorem establishes this conclusion in arbitrary dimension
and codimension for every nonpositive space form.

\begin{theorem}\label{thm:main}
	There are no proper biharmonic PNMCV submanifolds in
	$N^{m+p}(c)$ with~$c\leq0$.
\end{theorem}

For $c=0$, Theorem~\ref{thm:main} settles the PNMCV version of Chen's
conjecture formulated in~\cite{Chen2019}. For $c<0$, it establishes the generalized
Chen's conjecture for hyperbolic space under the PNMCV assumption.

For a surface in $N^4(c)$, the Ricci equation shows that the normal
connection is flat if and only if the shape operators commute, which does not imply that the mean curvature vector field  is
parallel. Very recently, Andronic et al.~\cite{AMOS2026} studied
biconservative surfaces with flat normal bundle under a Weingarten
relation between the mean curvature function $f$ and the Gaussian
curvature $K$. The authors derived an ordinary differential system for a class
of non-CMC biconservative surfaces and excluded biharmonicity in several
cases. A further classification theorem in~\cite{AMOS2026} depends on
a nonexistence conjecture for non-CMC biharmonic Weingarten surfaces.
Our second result needs only flatness of the normal connection. We prove 

\begin{theorem}\label{thm:flat}
	Every biharmonic surface with flat normal bundle in $N^4(c)$
	with~$c\leq0$ is minimal.
\end{theorem}

\begin{remark}
	Theorem~\ref{thm:flat} holds in both $\mathbb R^4$ and
	$\mathbb H^4(c)$, without a Weingarten relation or a PNMCV assumption.
	The proof relies on both the tangent and normal bundles having rank
	two. The results give important support for Chen's conjecture and the generalized Chen's conjecture  in codimension two. Note that the corresponding problem for biharmonic submanifolds in higher dimensions or higher codimensions remains a nontrivial open question, requiring further investigation.
\end{remark}

For a biharmonic immersion with constant mean curvature, minimality
follows directly from the normal biharmonic equation. Taking the inner
product with $H$, we obtain
\[
 \frac12\operatorname{div}\left(\nabla f^2\right)
 =\left|\nabla^\perp H\right|^2+\left|A_H\right|^2-mcf^2
 \geq mf^2\left(f^2-c\right).
\]
Since $c\leq0$, the function $f$ cannot be a positive constant.
It therefore remains to exclude a connected open set on which
$f>0$ and $\nabla f\neq0$.

For PNMCV submanifolds, we extend the continuation argument
in~\cite{FHY2026} to higher codimension. The tangential biharmonic
equation shows that $\nabla f$ is a principal direction of $A_\xi$.
The Gauss, Codazzi and Ricci equations then yield a differential system
along the integral curves of $\nabla f/|\nabla f|$.
The coefficients of $h$ in the normal directions orthogonal to $\xi$
satisfy quadratic first integrals involving the principal curvatures
of $A_\xi$ and the connection coefficients.
These first integrals and the trace relation give bounds on the
maximal interval of existence. For $c<0$,
Lemma~\ref{lem:growth-comparison} yields a contradiction.

For surfaces with flat normal bundle, we first prove that
$\nabla f$ is a common principal direction of the shape operators.
In a common principal frame, the Gauss, Codazzi and biharmonic
equations, together with the commutator equations, imply that one
principal derivative of $f$ vanishes.
The integral curves of $\nabla f/|\nabla f|$ are then geodesics, and
the curvature coefficients depend on one parameter.
Choosing a parallel orthonormal normal frame gives an autonomous
system of scalar equations for the mean curvature components,
the coefficients of $h$ and a connection coefficient.
Two first integrals and a monotone scalar quantity give a contradiction
on the maximal interval of existence.

The paper is organized as follows.
In Section~\ref{sec:preliminaries}, we recall the fundamental equations
and prove the integral estimate used in the continuation arguments.
In Section~\ref{sec:pnmc}, we prove Theorem~\ref{thm:main}.
In Sections~\ref{sec:flat-reduction} and~\ref{sec:vector}, we establish
the principal-direction reduction and complete the proof of
Theorem~\ref{thm:flat}.

\section{Preliminaries}\label{sec:preliminaries}

Let $\varphi\colon (M^m,g)\to N^{m+p}(c)$ be an isometric immersion
of codimension $p$ into a space form of constant sectional curvature $c$.
We denote by $\nabla$ and $\widetilde\nabla$ the Levi--Civita connections
of $M^m$ and $N^{m+p}(c)$, respectively. The second fundamental form is
denoted by $h$. The Gauss and Weingarten formulas are given by
\begin{equation*}
\begin{aligned}
 \widetilde\nabla_XY&=\nabla_XY+h(X,Y),\\
 \widetilde\nabla_X\eta&=-A_\eta X+\nabla^\perp_X\eta,
\end{aligned}
\end{equation*}
where $\nabla^\perp$ denotes the normal connection and $A_\eta$ denotes
the shape operator in the normal direction $\eta$.

Choose local orthonormal tangent and normal frames
$\{e_1,\ldots,e_m\}$ and $\{\nu_1,\ldots,\nu_p\}$.
Unless otherwise stated, the indices $i,j,k,l$ range from $1$ to $m$,
and $\alpha,\beta$ range from $1$ to $p$.
The components of $h$ and the mean curvature vector field $H$ are given by
\[
 h_{ij}^{\alpha}
 =\left\langle h\left(e_i,e_j\right),\nu_\alpha\right\rangle
 =\left\langle A_{\nu_\alpha}e_i,e_j\right\rangle,
 \quad H_\alpha=\frac1m\sum_{i=1}^m h_{ii}^{\alpha}.
\]
Thus
\[
 H=\sum_{\alpha=1}^p H_\alpha\nu_\alpha,
 \quad |H|^2=\sum_{\alpha=1}^p H_\alpha^2.
\]
The tangent and normal connection forms are defined by
\begin{align*}
 \omega_{i}^{j}\left(X\right)
 =\left\langle\nabla_Xe_i,e_j\right\rangle,\quad
 \omega_{m+\alpha}^{m+\beta}\left(X\right)
 =\left\langle\nabla_X^\perp\nu_\alpha,\nu_\beta\right\rangle.
\end{align*}
Since the frames are orthonormal, both families are skew-symmetric
in the two indices. The covariant derivatives of $h_{ij}^{\alpha}$
are given by
\begin{align*}
 h_{ijk}^{\alpha}
 ={}&e_k\left(h_{ij}^{\alpha}\right)
 +\sum_{\beta=1}^p h_{ij}^{\beta}
       \omega_{m+\beta}^{m+\alpha}\left(e_k\right)
 -\sum_{l=1}^m h_{lj}^{\alpha}\omega_{i}^{l}\left(e_k\right)
  -\sum_{l=1}^m h_{il}^{\alpha}\omega_{j}^{l}\left(e_k\right).
\end{align*}

We use the curvature convention
\[
 R(X,Y)Z=\nabla_X\nabla_YZ-\nabla_Y\nabla_XZ-\nabla_{[X,Y]}Z.
\]
Writing $R_{ijkl}=\left\langle R\left(e_i,e_j\right)e_k,e_l\right\rangle$,
the Gauss, Codazzi and Ricci equations are given by
\begin{align}
 &R_{ijkl}
 =c\left(\delta_{jk}\delta_{il}-\delta_{ik}\delta_{jl}\right)
   +\sum_{\alpha=1}^p
       \left(h_{il}^{\alpha}h_{jk}^{\alpha}
             -h_{ik}^{\alpha}h_{jl}^{\alpha}\right),
       \label{eq:Gauss-general}\\
 &\quad h_{ijk}^{\alpha}=h_{ikj}^{\alpha},
       \label{eq:Codazzi-general}\\
 &\left\langle R^\perp\left(e_i,e_j\right)\nu_\alpha,
                  \nu_\beta\right\rangle
 =\sum_{k=1}^m\left(h_{ik}^{\beta}h_{jk}^{\alpha}
                       -h_{ik}^{\alpha}h_{jk}^{\beta}\right).
       \label{eq:Ricci-general}
\end{align}
For a normal vector field $\eta$, the Codazzi equation can also be written as
\begin{equation}\label{eq:Codazzi-shape}
 \left(\nabla_XA_\eta\right)Y-\left(\nabla_YA_\eta\right)X
 =A_{\nabla^\perp_X\eta}Y-A_{\nabla^\perp_Y\eta}X.
\end{equation}
Equation~\eqref{eq:Ricci-general} implies that the normal connection
is flat if and only if the shape operators commute at every point.

For a smooth function $\psi$ on $M^m$,  the Laplacian acting on $\psi$ is defined by
\[
 \Delta \psi=\sum_{i=1}^m
 \left(e_i\left(e_i(\psi)\right)
       -\left(\nabla_{e_i}e_i\right)(\psi)\right).
\]
We use the same trace convention for the normal Laplacian.
For the mean curvature components, put
\begin{align*}
 H_{\alpha,i}
 &=e_i\left(H_\alpha\right)
   +\sum_{\beta=1}^p H_\beta
       \omega_{m+\beta}^{m+\alpha}\left(e_i\right),\\
 H_{\alpha,ij}
 &=e_j\left(H_{\alpha,i}\right)
   +\sum_{\beta=1}^p H_{\beta,i}
       \omega_{m+\beta}^{m+\alpha}\left(e_j\right)
   -\sum_{k=1}^m H_{\alpha,k}\omega_{i}^{k}\left(e_j\right).
\end{align*}
The normal and tangential biharmonic equations have the following
component form~(cf. \cite{OuChen2020,FetcuOniciuc2022}).

\begin{proposition}
The immersion $\varphi\colon M^m\to N^{m+p}(c)$ is biharmonic if and only if
\begin{align}
 &\sum_{i=1}^m H_{\alpha,ii}
 -\sum_{i,j=1}^m\sum_{\beta=1}^p
       H_\beta h_{ij}^{\beta}h_{ij}^{\alpha}
 +mcH_\alpha=0,
 \qquad 1\leq\alpha\leq p,
 \label{eq:biharmonic-normal-general}\\
 &2\sum_{i=1}^m\sum_{\alpha=1}^p
       h_{ij}^{\alpha}H_{\alpha,i}
 +\frac m2e_j\left(\sum_{\alpha=1}^pH_\alpha^2\right)=0,
 \qquad 1\leq j\leq m.
 \label{eq:biharmonic-tangent-general}
\end{align}
\end{proposition}

Define
\begin{equation}\label{eq:VU-sets}
 \begin{gathered}
  \mathcal V=\{x\in M:H(x)\ne0\},\qquad
  \mathcal U=\{x\in\mathcal V:\nabla f(x)\ne0\},\\
  f=|H|,\qquad \xi=\frac{H}{f}.
 \end{gathered}
\end{equation}
The function $f$ is smooth on $\mathcal V$.
Choose $\nu_1=\xi$ locally on $\mathcal V$, so that $H_1=f$ and
$H_\alpha=0$ for $2\leq\alpha\leq p$.
Setting $\alpha=1$ in \eqref{eq:biharmonic-normal-general} gives
\begin{equation}\label{eq:radial-general}
 \Delta f=f\left(
       \sum_{i,j=1}^m\left(h_{ij}^{1}\right)^2
       +\sum_{i=1}^m\sum_{\alpha=2}^p
          \left(\omega_{m+1}^{m+\alpha}\left(e_i\right)\right)^2
       -mc\right).
\end{equation}
Since $\sum_i h_{ii}^{1}=mf$, the Cauchy--Schwarz inequality yields
\[
 |A_\xi|^2=\sum_{i,j=1}^m\left(h_{ij}^{1}\right)^2\geq mf^2.
\]
For $c\leq0$, equation~\eqref{eq:radial-general} therefore excludes
a positive constant value of $f$ on a nonempty open set.

On an open subset of $\mathcal V$ where $\nabla^\perp\xi=0$,
equation~\eqref{eq:Codazzi-shape} implies that $A_\xi$ is a Codazzi
tensor. Since $R^\perp(X,Y)\xi=0$, equation~\eqref{eq:Ricci-general}
also gives
\[
 [A_\xi,A_\eta]=0
 \qquad\text{for every normal field }\eta.
\]

We use the following estimate in the continuation arguments in
Sections~\ref{sec:pnmc} and~\ref{sec:vector}. All derivatives in the
statement are taken with respect to $t$.

\begin{lemma}\label{lem:growth-comparison}
Let $w\in C^1([0,\infty))$ be positive, and suppose that
\[
 w(t)\leq C_0(1+t)^{r_0} e^{at},\qquad t\geq0,
\]
for constants $C_0>0$ and $a,r_0\geq0$. There is no bounded function
$\zeta\in C^2([0,\infty))$ satisfying
\[
 \zeta'>0,\qquad \left(\zeta'w\right)'\geq b_0 w
\]
for a constant $b_0>0$.
\end{lemma}

\begin{proof}
Suppose that such a function exists and put
$W(t)=\int_0^t w(s)\,\mathrm ds$. Integrating the differential
inequality for $\zeta'w$ and using $\zeta'(0)w(0)>0$, we obtain
\[
 \zeta'(t)\geq b_0\frac{W(t)}{w(t)}.
\]
Since $\zeta$ is increasing and bounded, it follows that
\[
 L:=\int_1^\infty\frac{W(t)}{w(t)}\,\mathrm dt<\infty.
\]
For $T>1$, the Cauchy--Schwarz inequality yields
\begin{equation}\label{eq:growth-CS}
\begin{aligned}
 (T-1)^2
 \leq\left(\int_1^T\frac{W}{w}\,\mathrm dt\right)
       \left(\int_1^T\frac{w}{W}\,\mathrm dt\right)
 \leq L\log\frac{W(T)}{W(1)}.
\end{aligned}
\end{equation}
The bound for $w$ gives
$W(T)\leq C_1(1+T)^{r_0+1}e^{aT}$ for some $C_1>0$.
The right-hand side of \eqref{eq:growth-CS} is therefore
$O(T+\log(1+T))$, contradicting the quadratic growth of the
left-hand side as $T\to\infty$.
\end{proof}

\section{PNMCV submanifolds}\label{sec:pnmc}

Let $\varphi\colon M^m\to N^{m+p}(c)$ be a biharmonic PNMCV
immersion with $c\leq0$. Since $f>0$ and $\nabla^\perp\xi=0$
on $M$, we have $\mathcal V=M$ in \eqref{eq:VU-sets}.
Put $A=A_\xi$. Choose local orthonormal tangent and normal frames
$\{e_1,\ldots,e_m\}$ and $\{\nu_1=\xi,\nu_2,\ldots,\nu_p\}$.
The components of the second fundamental form are given by
$h_{ij}^{\alpha}=\langle h(e_i,e_j),\nu_\alpha\rangle$.
The component of \eqref{eq:biharmonic-normal-general} in the
direction $\xi$ and the tangential equation
\eqref{eq:biharmonic-tangent-general} give
\begin{align}
&\Delta f=f\left(\sum_{i,j=1}^m\left(h_{ij}^{1}\right)^2-mc\right),
\label{eq:biharmonic}\\
&2\sum_{j=1}^m h_{ij}^{1}e_j(f)+mf e_i(f)=0,
\qquad 1\leq i\leq m.
\label{eq:biharmonic-tangent}
\end{align}
The remaining normal equations are given by
\begin{equation}\label{eq:mixed-normal}
\sum_{i,j=1}^m h_{ij}^{1}h_{ij}^{\alpha}=0,
\qquad 2\leq\alpha\leq p.
\end{equation}
Since $\nabla^\perp\xi=0$, equations
\eqref{eq:Codazzi-shape} and \eqref{eq:Ricci-general} imply
\begin{equation}\label{eq:codazzi-A}
\left(\nabla_XA\right)Y=\left(\nabla_YA\right)X
\end{equation}
and
\begin{equation}\label{eq:commuting}
[A,A_\eta]=0
\qquad\text{for every normal field }\eta,
\end{equation}
since $R^\perp(X,Y)\xi=0$.

Assume $m\geq2$ and $\mathcal U\ne\varnothing$. The number of distinct
eigenvalues of $A$ takes values in $\{1,\ldots,m\}$. Choose a point
of $\mathcal U$ at which this number is maximal on $\mathcal U$.
The eigenvalue gaps at that point persist on a sufficiently small
connected neighborhood contained in $\mathcal U$. Maximality shows that
none of the resulting eigenvalue groups can split on this neighborhood.
Thus the number of distinct eigenvalues and the multiplicities are
constant there. If an eigenvalue has multiplicity $m_0$, the implicit
function theorem applied to the $(m_0-1)$-st derivative of the
characteristic polynomial gives a smooth local eigenvalue function.
The corresponding eigenspaces form smooth distributions.
For the local description of principal curvatures and the corresponding
eigendistributions, see~\cite{KN1969}.
On this neighborhood, set
\begin{equation}\label{eq:distinguished}
e_1=\frac{\nabla f}{|\nabla f|},\qquad
\kappa_1=-\frac m2f.
\end{equation}
By \eqref{eq:biharmonic-tangent}, $Ae_1=\kappa_1 e_1$. Moreover,
$e_1(\kappa_1)\ne0$ and $X(\kappa_1)=0$ for $X\perp e_1$.
For the remainder of this section, indices $i,j$ range from $2$ to $m$
and the normal index $\alpha$ ranges from $2$ to $p$, unless a different
range is specified. A sum or product over an empty index set is understood
to be $0$ or $1$, respectively.

\begin{lemma}\label{lem:frame}
The eigenvalue $\kappa_1$ is simple and $\nabla_{e_1}e_1=0$.
Let $\mathcal D=e_1^\perp$ and define
\[
\mathscr{S}X=-\nabla_Xe_1,\qquad X\in\mathcal D.
\]
Then $\mathscr{S}$ is self-adjoint, $[A|_{\mathcal D},\mathscr{S}]=0$, and
\begin{equation}\label{eq:derivative-A}
\left(\nabla_{e_1}A\right)X=(A-\kappa_1 \mathrm{Id})\mathscr{S}X,
\qquad X\in\mathcal D,
\end{equation}
where $\mathrm{Id}$ denotes the identity endomorphism.
Along an integral curve of $e_1$, there is an orthonormal frame
$e_2,\ldots,e_m$ of $\mathcal D$ such that
\begin{equation}\label{eq:parallel-frame}
\nabla_{e_1}e_i=0,\qquad
Ae_i=\kappa_i e_i,\qquad \mathscr{S}e_i=b_i e_i,
\end{equation}
where
$
 b_i=\left\langle\nabla_{e_i}e_i,e_1\right\rangle
 =\omega_{i}^{1}\left(e_i\right).
$
Thus the functions $b_i$ are the principal curvatures of the level hypersurfaces
of $f$ with respect to the unit normal $e_1$. In addition,
\begin{equation}\label{eq:mixed-B}
h_{1i}^{\beta}=0,\qquad 2\leq i\leq m,\quad 1\leq\beta\leq p.
\end{equation}
\end{lemma}

\begin{proof}
For a Codazzi tensor, an eigenvalue of multiplicity at least two is
constant in the directions of the corresponding eigenspace. Indeed, let
$X,Y$ be orthonormal local fields in the same eigenspace, with eigenvalue
$\lambda_1$. Equation~\eqref{eq:codazzi-A} gives
\[
X\left(\lambda_1\right)
=\left\langle\left(\nabla_XA\right)Y,Y\right\rangle
=\left\langle\left(\nabla_YA\right)X,Y\right\rangle=0.
\]
Since $e_1(\kappa_1)\ne0$, the eigenvalue $\kappa_1$ is simple.

For a local unit eigenvector field $e_i$ orthogonal to $e_1$
and corresponding to the eigenvalue $\kappa_i$, \eqref{eq:codazzi-A} yields
\[
(\kappa_1-\kappa_i)
\left\langle\nabla_{e_1}e_1,e_i\right\rangle
=\left\langle\left(\nabla_{e_i}A\right)e_1,e_1\right\rangle
=e_i(\kappa_1)=0.
\]
Thus $\nabla_{e_1}e_1=0$. The distribution $\mathcal D$ is tangent
to the level sets of $\kappa_1$, so $\mathscr{S}$ is self-adjoint.

For $X,Y\in\mathcal D$, taking the $e_1$ component in
\eqref{eq:codazzi-A} gives
\[
\langle(A-\kappa_1 \mathrm{Id})Y,\mathscr{S}X\rangle
=\langle(A-\kappa_1 \mathrm{Id})X,\mathscr{S}Y\rangle.
\]
Since $A$ and $\mathscr{S}$ are self-adjoint, it follows that they commute
on $\mathcal D$. Applying \eqref{eq:codazzi-A} to $e_1$ and $X$ gives
\[
\left(\nabla_{e_1}A\right)X
=\left(\nabla_XA\right)e_1
=(A-\kappa_1 \mathrm{Id})\mathscr{S}X,
\]
which yields \eqref{eq:derivative-A}.

The right-hand side of \eqref{eq:derivative-A} commutes with
$A|_{\mathcal D}$. Let $X$ and $Y$ belong to eigenspaces with distinct
eigenvalues $\lambda_1$ and $\lambda_2$. Differentiating $AX=\lambda_1 X$
along $e_1$ and using \eqref{eq:derivative-A}, we obtain
\[
(\lambda_1-\lambda_2)\left\langle\nabla_{e_1}X,Y\right\rangle
=\left\langle\left(\nabla_{e_1}A\right)X,Y\right\rangle=0.
\]
Since $\nabla_{e_1}e_1=0$, parallel transport along $e_1$ also preserves
$\mathcal D$. Hence the eigenspaces of $A|_{\mathcal D}$ are preserved
by parallel transport.
On an eigenspace with eigenvalue $\lambda_1$, equation
\eqref{eq:derivative-A} becomes
\[
e_1\left(\lambda_1\right)\mathrm{Id}=(\lambda_1-\kappa_1)\mathscr{S}.
\]
Since $\lambda_1\ne\kappa_1$, it follows that $\mathscr{S}$ is a scalar multiple of
the identity on each such eigenspace. Parallel transport of an orthonormal
eigenbasis along an integral curve gives \eqref{eq:parallel-frame}.

Finally, \eqref{eq:commuting} and the simplicity of $\kappa_1$
show that $e_1$ is an eigenvector of every $A_\eta$. Hence
\[
h_{1i}^{\beta}=\left\langle A_{\nu_\beta}e_1,e_i\right\rangle=0,
\qquad 1\leq\beta\leq p,
\]
which gives \eqref{eq:mixed-B}.
\end{proof}

Fix an integral curve of $e_1$, parametrized by arc length $t$,
and use the frame defined in Lemma~\ref{lem:frame}. Choose the orthonormal normal
frame $\{\nu_1=\xi,\nu_2,\ldots,\nu_p\}$ to be parallel along this
curve. A prime denotes differentiation with respect to $t$.
In this frame, we have
\begin{equation}\label{eq:normal-vectors}
h_{rs}^{1}=\kappa_r\delta_{rs},\qquad 1\leq r,s\leq m.
\end{equation}
Thus $\kappa_r=h_{rr}^{1}$ for $1\leq r\leq m$. For $\alpha\geq2$,
the components $h_{ij}^{\alpha}$ with $i,j\geq2$ may be nonzero even
when $i\ne j$.
Since $H=f\xi$, we have $\tr A_\xi=mf$ and
$\tr A_{\nu_\alpha}=0$ for $\alpha\geq2$. Using
$\kappa_1=-mf/2$ from \eqref{eq:distinguished}, we obtain
\begin{equation*}
\sum_{i=2}^m\kappa_i=-3\kappa_1,
\qquad
\sum_{i=2}^m h_{ii}^{\alpha}=-h_{11}^{\alpha},
\qquad 2\leq\alpha\leq p.
\end{equation*}

\begin{proposition}\label{prop:reduction}
Along the integral curve of $e_1$, we have
\begin{align}
&\kappa_i'=(\kappa_i-\kappa_1)b_i,
\label{eq:vector-system}\\
&\left(h_{ii}^{\alpha}\right)'=
\left(h_{ii}^{\alpha}-h_{11}^{\alpha}\right)b_i,
\label{eq:normal-component-system}\\
&b_i'=b_i^2+\kappa_1\kappa_i
+\sum_{\alpha=2}^p h_{11}^{\alpha}h_{ii}^{\alpha}+c,
\qquad 2\leq i\leq m.
\label{eq:connection-system}
\end{align}
Moreover,
\begin{equation}\label{eq:scalar-system}
\kappa_1''=\left(\sum_{i=2}^m b_i\right)\kappa_1'
+\kappa_1\left(\sum_{j=1}^m\kappa_j^2-mc\right).
\end{equation}
\end{proposition}

\begin{proof}
Applying \eqref{eq:derivative-A} to $e_i$ and taking the $e_i$
component gives \eqref{eq:vector-system}. By \eqref{eq:parallel-frame} and
\eqref{eq:mixed-B}, we have $\nabla_{e_1}e_i=0$,
$\nabla_{e_i}e_1=-b_i e_i$, and $h_{1i}^{\beta}=0$ for every
normal index $\beta$. Since the normal frame is parallel along the
curve, the normal connection coefficients in the direction $e_1$
vanish. Substituting into \eqref{eq:Codazzi-general}, we obtain
\[
\left(h_{ii}^{\alpha}\right)'
=b_i h_{ii}^{\alpha}-b_i h_{11}^{\alpha},
\]
which gives \eqref{eq:normal-component-system}.

Since $\nabla_{e_1}e_1=0$ and $\nabla_{e_1}e_i=0$, the sectional
curvature of the plane spanned by $e_i,e_1$ is given by
\[
K(e_i,e_1)=b_i'-b_i^2.
\]
On the other hand, combining \eqref{eq:Gauss-general} with \eqref{eq:mixed-B} yields
\[
\begin{aligned}
K(e_i,e_1)
&=c+\sum_{\beta=1}^p\left(h_{11}^{\beta}h_{ii}^{\beta}
-\left(h_{1i}^{\beta}\right)^2\right)\\
&=c+\kappa_1\kappa_i
+\sum_{\alpha=2}^p h_{11}^{\alpha}h_{ii}^{\alpha}.
\end{aligned}
\]
Comparing the two expressions for $K(e_i,e_1)$ gives
\eqref{eq:connection-system}.

Since $e_i(f)=0$ for $2\leq i\leq m$ and
$\left\langle\nabla_{e_i}e_i,e_1\right\rangle=b_i$, we have
\[
\Delta f=f''-\left(\sum_{i=2}^m b_i\right)f'.
\]
Substituting this expression for $\Delta f$ and
\eqref{eq:normal-vectors} into \eqref{eq:biharmonic}, and using
$\kappa_1=-mf/2$, we obtain \eqref{eq:scalar-system}.
\end{proof}

Let $m\geq2$, $p\geq1$, and $c\leq0$. We consider
\eqref{eq:vector-system}--\eqref{eq:connection-system} as a system of ordinary differential equations
for the real functions $\kappa_i$, $b_i$, and $h_{ii}^{\alpha}$, where
\begin{equation}\label{eq:elimination-traces}
\kappa_1=-\frac13\sum_{i=2}^m\kappa_i,
\qquad
h_{11}^{\alpha}=-\sum_{i=2}^m h_{ii}^{\alpha},
\qquad 2\leq\alpha\leq p.
\end{equation}
Substituting \eqref{eq:elimination-traces} into
\eqref{eq:vector-system}--\eqref{eq:connection-system} gives a
polynomial system in these unknown functions.
Thus every local solution extends uniquely to a maximal open interval
$I$, and the extended solution is real analytic. The trace equations
\eqref{eq:elimination-traces} hold throughout $I$ by definition.
Moreover, $\kappa_1$ and all derivatives of $\kappa_1$ are analytic on
$I$. If \eqref{eq:scalar-system} holds on a nonempty open subinterval,
the difference between the two sides is analytic on $I$ and vanishes on
that subinterval. The identity theorem therefore gives
\eqref{eq:scalar-system} on all of $I$. We restrict to solutions
satisfying \eqref{eq:scalar-system} on $I$.

For the local solution obtained in Proposition~\ref{prop:reduction},
equation~\eqref{eq:scalar-system} therefore holds on the maximal
interval $I$.

After translating the parameter, we may assume that $0\in I$.
Following the argument in~\cite{FHY2026}, put $c_0=\sqrt{-c}$ and
define $\rho_i$ by
\begin{equation}\label{eq:rho-initial}
\begin{aligned}
\rho_i'&=-b_i\rho_i,\\
\rho_i(0)&=
\left(1+b_i(0)^2+\kappa_i(0)^2
+\sum_{\alpha=2}^p\left(h_{ii}^{\alpha}(0)\right)^2+c_0^2\right)^{-1/2}.
\end{aligned}
\end{equation}
Since the first-order equations in \eqref{eq:rho-initial} are linear
and the initial values are positive, the functions $\rho_i$ are
strictly positive on $I$.

\begin{lemma}\label{lem:identities}
The following equations hold on $I$:
\begin{align}
&\left[\left(b_i^2+\kappa_i^2+\sum_{\alpha=2}^p\left(h_{ii}^{\alpha}\right)^2+c\right)\rho_i^2\right]'
=0,\label{eq:energy}\\
&\left[(\kappa_1-\kappa_i)\rho_i\right]'
=\kappa_1'\rho_i,\label{eq:product-single}\\
&\left(\kappa_1'\prod_{i=2}^m\rho_i\right)'
=\kappa_1\left(\sum_{j=1}^m\kappa_j^2-mc\right)\prod_{i=2}^m\rho_i.
\label{eq:product-all}
\end{align}
Moreover,
\begin{equation}\label{eq:energy-bound}
\left(b_i^2+\kappa_i^2+\sum_{\alpha=2}^p\left(h_{ii}^{\alpha}\right)^2\right)\rho_i^2
\leq1+c_0^2\rho_i^2,
\end{equation}
and consequently
\begin{equation}\label{eq:linear-bound}
|b_i|\rho_i\leq1+c_0\rho_i,
\qquad
|\kappa_i|\rho_i\leq1+c_0\rho_i,
\qquad
\left|h_{ii}^{\alpha}\right|\rho_i\leq1+c_0\rho_i.
\end{equation}
\end{lemma}

\begin{proof}
Differentiating and using
\eqref{eq:vector-system}--\eqref{eq:connection-system}, we obtain
\begin{align*}
&\left(b_i^2+\kappa_i^2
+\sum_{\alpha=2}^p\left(h_{ii}^{\alpha}\right)^2+c\right)'\\
&\quad=2b_i\left(b_i^2+\kappa_1\kappa_i
+\sum_{\alpha=2}^p h_{11}^{\alpha}h_{ii}^{\alpha}+c\right)
+2\kappa_i(\kappa_i-\kappa_1)b_i\\
&\qquad\,+2b_i\sum_{\alpha=2}^p
h_{ii}^{\alpha}\left(h_{ii}^{\alpha}-h_{11}^{\alpha}\right)\\
&\quad=2b_i\left(b_i^2+\kappa_i^2
+\sum_{\alpha=2}^p\left(h_{ii}^{\alpha}\right)^2+c\right).
\end{align*}
Since $\rho_i'=-b_i\rho_i$, equation~\eqref{eq:energy} follows.
Evaluating the constant expression in \eqref{eq:energy} at $t=0$
and using \eqref{eq:rho-initial}, we find that this constant is at
most $1$. Hence \eqref{eq:energy-bound} follows.
Since $\sqrt{1+c_0^2\rho_i^2}\leq1+c_0\rho_i$,
equation~\eqref{eq:energy-bound} gives \eqref{eq:linear-bound}.

Next, using \eqref{eq:vector-system} and \eqref{eq:rho-initial}, we obtain
\begin{align*}
\left[(\kappa_1-\kappa_i)\rho_i\right]'
&=\left(\kappa_1'-(\kappa_i-\kappa_1)b_i\right)\rho_i
-(\kappa_1-\kappa_i)b_i\rho_i\\
&=\kappa_1'\rho_i.
\end{align*}
Differentiating the product and using \eqref{eq:rho-initial} yields
\[
\left(\prod_i\rho_i\right)'=-\left(\sum_{i=2}^m b_i\right)\prod_i\rho_i,
\]
so \eqref{eq:scalar-system} implies \eqref{eq:product-all}.
\end{proof}

\begin{lemma}\label{lem:bounded-factors}
Let $I$ be the maximal interval of a solution of
\eqref{eq:vector-system}--\eqref{eq:connection-system} with \eqref{eq:elimination-traces}, and
suppose that \eqref{eq:scalar-system} holds on $I$. Put
$t_*=\sup I$. There is no $t_0\in I$ such that, on
$[t_0,t_*)$,
\begin{equation}\label{eq:positive-growth}
\kappa_1>0,\qquad \kappa_1'>0,\qquad
\kappa_1'\prod_i\rho_i\geq\delta>0,
\end{equation}
and
\begin{equation}\label{eq:rho-upper}
0<\rho_i\leq C_i,\qquad 2\leq i\leq m,
\end{equation}
where $\delta,C_2,\ldots,C_m$ are constants.
\end{lemma}

\begin{proof}
Suppose that \eqref{eq:positive-growth} and \eqref{eq:rho-upper}
hold. By \eqref{eq:energy-bound}, there is a constant $R_0>0$ such
that
\begin{equation}\label{eq:weighted-bound}
|b_i|\rho_i\leq R_0,\qquad
|\kappa_i|\rho_i\leq R_0,\qquad
\left|h_{ii}^{\alpha}\right|\rho_i\leq R_0.
\end{equation}
For each $t\in[t_0,t_*)$, equation~\eqref{eq:elimination-traces}
implies that there is an index $i$ such that
\[
\kappa_i(t)\leq-\frac{3\kappa_1(t)}{m-1}<0.
\]
For this index, \eqref{eq:weighted-bound} gives
\begin{equation}\label{eq:fixed-index-bound}
0<(\kappa_1-\kappa_i)\rho_i
=\left(1+\frac{\kappa_1}{|\kappa_i|}\right)|\kappa_i|\rho_i
\leq\frac{m+2}{3}R_0=:K_0.
\end{equation}

Choose a strictly increasing sequence in $[t_0,t_*)$ tending to $t_*$.
Since there are only finitely many indices, one fixed index satisfies
\eqref{eq:fixed-index-bound} along an infinite subsequence, which still
tends to $t_*$. Relabel this index as $2$ and put
\[
q=(\kappa_1-\kappa_2)\rho_2.
\]
By \eqref{eq:product-single}, we have $q'=\kappa_1'\rho_2>0$.
For each $t\in[t_0,t_*)$, choose a point of the subsequence beyond
$t$. Since $q$ is increasing, it follows that $q(t)\leq K_0$.
Thus $q\leq K_0$ throughout $[t_0,t_*)$. Taking $t_1$ from the
subsequence and setting $\varepsilon=q(t_1)>0$, we obtain
\begin{equation}\label{eq:v-two-sided}
\varepsilon\leq q(t)\leq K_0,\qquad t_1\leq t<t_*.
\end{equation}

For $j\ne2$, equations \eqref{eq:positive-growth},
\eqref{eq:rho-upper}, and \eqref{eq:weighted-bound} give
\begin{align*}
\left|(\log\rho_j)'\right|
&=|b_j|\leq\frac{R_0}{\rho_j}
\leq\frac{R_0}{\delta}\kappa_1'\prod_{i\ne j}\rho_i\\
&\leq\frac{R_0}{\delta}
\left(\prod_{i\ne2,j}C_i\right)\kappa_1'\rho_2
=D_jq',
\end{align*}
where $D_j$ is a fixed positive constant. Integrating from $t_1$ to
$t$ and using \eqref{eq:v-two-sided}, we obtain
\begin{equation}\label{eq:rho-other-lower}
\left|\log\frac{\rho_j(t)}{\rho_j(t_1)}\right|
\leq D_j(K_0-\varepsilon),
\qquad t_1\leq t<t_*.
\end{equation}
Thus each $\rho_j$, $j\ne2$, has a uniform positive lower bound.
It follows from \eqref{eq:weighted-bound} that
$\kappa_j$, $b_j$, and $h_{jj}^{\alpha}$ are uniformly bounded for $j\ne2$.

Choose $C_0\geq0$ such that
$
\left|\sum_{j\ne2}\kappa_j\right|\leq C_0.
$
Using \eqref{eq:elimination-traces} again, we obtain
\begin{equation}\label{eq:lambda-difference}
0<\kappa_1-\kappa_2
=4\kappa_1+\sum_{j\ne2}\kappa_j
\leq4\kappa_1+C_0.
\end{equation}
Since both denominators below are positive,
\eqref{eq:lambda-difference} and \eqref{eq:product-single} give
\[
\frac{\kappa_1'}{4\kappa_1+C_0}
\leq\frac{\kappa_1'}{\kappa_1-\kappa_2}
=\frac{q'}q.
\]
Integrating from $t_1$ to $t$ and using \eqref{eq:v-two-sided}, we obtain
\begin{equation}\label{eq:lambda-log-bound}
\frac14\log\frac{4\kappa_1(t)+C_0}{4\kappa_1(t_1)+C_0}
\leq\log\frac{q(t)}{\varepsilon}
\leq\log\frac{K_0}{\varepsilon}.
\end{equation}
It follows from \eqref{eq:lambda-log-bound} that $\kappa_1$ has a
uniform upper bound. Combining \eqref{eq:v-two-sided} with
\eqref{eq:lambda-difference}, we obtain
\[
\rho_2=
\frac{q}{\kappa_1-\kappa_2}
\geq\frac{\varepsilon}{4\kappa_1+C_0}.
\]
Since $\kappa_1$ is uniformly bounded, the above gives a uniform positive
lower bound for $\rho_2$. It follows from \eqref{eq:weighted-bound} that
$\kappa_i$, $b_i$, and $h_{ii}^{\alpha}$ are uniformly bounded on $[t_1,t_*)$.

If $t_*<\infty$, the right-hand sides of \eqref{eq:vector-system}--\eqref{eq:connection-system}
are bounded. Hence $\kappa_i$, $b_i$, and $h_{ii}^{\alpha}$ have bounded derivatives
and finite limits as $t\uparrow t_*$. Local existence with these
limiting values, followed by uniqueness, extends the solution
beyond $t_*$. This contradicts maximality. Thus $t_*=\infty$.

On the other hand, \eqref{eq:positive-growth} and
\eqref{eq:rho-upper} imply
\[
\kappa_1'\geq\frac{\delta}{\prod_iC_i}>0.
\]
Integrating this inequality shows that $\kappa_1$ is unbounded,
which contradicts \eqref{eq:lambda-log-bound}.
\end{proof}

\begin{proposition}\label{prop:ode-nonexistence}
Let $m\geq2$, $p\geq1$, and $c\leq0$. Every smooth real
solution of \eqref{eq:vector-system}--\eqref{eq:connection-system},
\eqref{eq:elimination-traces}, and \eqref{eq:scalar-system} on
a connected open interval satisfies $\kappa_1\equiv0$.
\end{proposition}

\begin{proof}
Suppose first that $\kappa_1$ is nonconstant. Choose a point at
which both $\kappa_1$ and $\kappa_1'$ are nonzero. Equations
\eqref{eq:vector-system}--\eqref{eq:connection-system},
\eqref{eq:elimination-traces}, and \eqref{eq:scalar-system} are
invariant under simultaneously replacing $\kappa_1,\kappa_i$
by $-\kappa_1,-\kappa_i$, with $b_i$ and $h_{ii}^{\alpha}$ unchanged.
These equations are also invariant under reversal of the parameter, accompanied
by changing the signs of the $b_i$. After applying these
transformations and translating the parameter, we may assume
\begin{equation}\label{eq:positive-initial}
\kappa_1(0)>0,\qquad\kappa_1'(0)>0.
\end{equation}

Extend the solution of
\eqref{eq:vector-system}--\eqref{eq:connection-system}, with
\eqref{eq:elimination-traces}, to the maximal interval $I$ and write
$t_*=\sup I$. Since the extended solution is analytic,
\eqref{eq:scalar-system} holds throughout $I$ by the identity theorem.
As long as $\kappa_1>0$, we have
$\sum_{j=1}^m\kappa_j^2-mc\geq\kappa_1^2>0$. Integrating
\eqref{eq:product-all} from $0$ to $t$ and using
\eqref{eq:positive-initial}, we obtain
\[
\kappa_1'(t)\prod_i\rho_i(t)
\geq\kappa_1'(0)\prod_i\rho_i(0)=:\delta>0.
\]
It follows that $\kappa_1'>0$ as long as $\kappa_1>0$.
Consequently, $\kappa_1$ cannot return to zero, and
\begin{equation}\label{eq:global-positive}
\kappa_1(t)\geq\kappa_1(0)>0,\qquad
\kappa_1'(t)>0,\qquad
\kappa_1'(t)\prod_i\rho_i(t)\geq\delta>0
\end{equation}
throughout $[0,t_*)$.

It follows from \eqref{eq:rho-initial} and \eqref{eq:linear-bound} that
\[
\rho_i'\leq1+c_0\rho_i.
\]
Multiplying by $e^{-c_0 t}$ and integrating from $0$ to $t$ yields
\begin{equation}\label{eq:rho-exponential}
\rho_i(t)\leq e^{c_0 t}\left(\rho_i(0)+t\right),
\qquad 0\leq t<t_*.
\end{equation}
If $t_*<\infty$, the functions $\rho_i$ are uniformly bounded above
on $[0,t_*)$. Together with \eqref{eq:global-positive}, this
contradicts Lemma~\ref{lem:bounded-factors}. Therefore $t_*=\infty$.

We also have
\begin{equation}\label{eq:lambda-curvature-bound}
0<\kappa_1(0)\leq\kappa_1(t)\leq c_0,
\qquad t\geq0.
\end{equation}
Suppose that $\kappa_1(t_0)>c_0$ at some $t_0\geq0$.
Since $\kappa_1'>0$, we then have $\kappa_1>c_0$ for all
$t\geq t_0$. Set
\[
q_i=(\kappa_1-\kappa_i)\rho_i.
\]
The estimate $\kappa_i\rho_i\leq1+c_0\rho_i$ in
\eqref{eq:linear-bound} gives
\begin{equation}\label{eq:shifted-positive}
(\kappa_1-c_0)\rho_i\leq q_i+1.
\end{equation}
Since $q_i+1>0$ on this interval and
$q_i'=\kappa_1'\rho_i$ by \eqref{eq:product-single}, differentiation
and \eqref{eq:linear-bound} give
\begin{align}
\left(\frac{q_i+1}{\kappa_1-c_0}\right)'
&=\frac{\kappa_1'}{(\kappa_1-c_0)^2}
\left((\kappa_i-c_0)\rho_i-1\right)
\leq0.
\label{eq:shifted-monotonicity}
\end{align}
Combining \eqref{eq:shifted-positive} and
\eqref{eq:shifted-monotonicity}, we obtain, for every $t\geq t_0$,
\begin{equation*}
0<\rho_i(t)
\leq\frac{q_i(t)+1}{\kappa_1(t)-c_0}
\leq\frac{q_i(t_0)+1}{\kappa_1(t_0)-c_0}.
\end{equation*}
Thus the functions $\rho_i$ are uniformly bounded above on
$[t_0,\infty)$. Together with \eqref{eq:global-positive}, this
contradicts Lemma~\ref{lem:bounded-factors} and gives
\eqref{eq:lambda-curvature-bound}. If $c=0$, then $c_0=0$,
and \eqref{eq:lambda-curvature-bound} is already impossible.

It remains to consider $c<0$. Put
\begin{equation*}
V(t)=\prod_{i=2}^m\rho_i(t)>0.
\end{equation*}
By \eqref{eq:rho-exponential}, there is a constant $C>0$ such that
\begin{equation}\label{eq:VW-growth}
V(t)\leq C(1+t)^{m-1}e^{(m-1)c_0 t},\qquad t\geq0.
\end{equation}
Moreover, \eqref{eq:product-all} and \eqref{eq:global-positive}
give
\begin{equation}\label{eq:pnmc-growth-comparison}
\left(\kappa_1'V\right)'=\kappa_1\left(\sum_{j=1}^m\kappa_j^2+mc_0^2\right)V
\geq mc_0^2\kappa_1(0)V.
\end{equation}
Since $mc_0^2\kappa_1(0)>0$, equations
\eqref{eq:VW-growth} and \eqref{eq:pnmc-growth-comparison}, together
with \eqref{eq:lambda-curvature-bound} and \eqref{eq:global-positive},
show that the hypotheses of Lemma~\ref{lem:growth-comparison} hold with
$w=V$, $\zeta=\kappa_1$, $r_0=m-1$, $a=(m-1)c_0$, and
$b_0=mc_0^2\kappa_1(0)$. The lemma gives a contradiction, so
$\kappa_1$ is constant.

Since $\kappa_1$ is constant, \eqref{eq:scalar-system} gives
$\kappa_1\left(\sum_{j=1}^m\kappa_j^2-mc\right)=0$. If $\kappa_1\ne0$, then
$\sum_{j=1}^m\kappa_j^2-mc\geq\kappa_1^2>0$, a contradiction. Therefore
$\kappa_1\equiv0$.
\end{proof}

\begin{proof}[Proof of Theorem~\ref{thm:main}]
Suppose that such an immersion exists. Since the mean curvature vector
field is nowhere zero, we have $f=|H|>0$ and $H=f\xi$ on $M$.
If $m=1$, then $h_{11}^{1}=f$, and the tangential equation
\eqref{eq:biharmonic-tangent} becomes $3f e_1(f)=0$. It follows that $f$
is locally constant. Substituting this into \eqref{eq:biharmonic}
gives $0=f\left(f^2-c\right)$, which contradicts $f>0$ and $c\leq0$.

Assume $m\geq2$ and suppose that $\mathcal U\ne\varnothing$,
where $\mathcal U$ is defined in \eqref{eq:VU-sets}.
On a sufficiently small connected open subset of $\mathcal U$,
Lemma~\ref{lem:frame} and Proposition~\ref{prop:reduction}
give an integral curve along which
\eqref{eq:vector-system}--\eqref{eq:connection-system},
\eqref{eq:elimination-traces}, and \eqref{eq:scalar-system} hold.
By \eqref{eq:distinguished}, we have
\[
 \kappa_1'=-\frac m2e_1(f)\ne0.
\]
On the other hand, Proposition~\ref{prop:ode-nonexistence} gives
$\kappa_1\equiv0$. This contradiction shows that
$\mathcal U=\varnothing$. Hence $\nabla f=0$ on $M$.

It follows that $f$ is constant on each connected component of $M$.
The normal equation~\eqref{eq:biharmonic} gives
\[
0=f\left(\sum_{i,j=1}^m\left(h_{ij}^{1}\right)^2-mc\right).
\]
Since $\sum_{i=1}^m h_{ii}^{1}=mf$, the Cauchy--Schwarz inequality yields
\[
\sum_{i,j=1}^m\left(h_{ij}^{1}\right)^2
\geq\sum_{i=1}^m\left(h_{ii}^{1}\right)^2
\geq\frac1m\left(\sum_{i=1}^m h_{ii}^{1}\right)^2
=mf^2.
\]
It follows that
\[
\sum_{i,j=1}^m\left(h_{ij}^{1}\right)^2-mc
\geq m\left(f^2-c\right)>0,
\]
which contradicts \eqref{eq:biharmonic} and $f>0$.
\end{proof}

\begin{remark}
Since $\nabla^\perp\xi=0$, equations \eqref{eq:codazzi-A} and
\eqref{eq:commuting} show that $A_\xi$ is a Codazzi tensor and
$[A_\xi,A_\eta]=0$ for every normal field $\eta$. Since $\kappa_1$
is simple, we obtain $h_{1i}^{\beta}=0$ for every normal index
$\beta$, as used in Proposition~\ref{prop:reduction}. The operators
$A_\eta$ with $\eta\perp\xi$ need not commute with one another.
\end{remark}

\section{Surfaces with flat normal bundle}\label{sec:flat-reduction}

Let $M^2$ be a biharmonic surface in $N^4(c)$, where $c\leq0$, with
flat normal bundle. On a sufficiently small connected open subset of
$\mathcal V=\{x\in M^2:H(x)\neq0\}$ from \eqref{eq:VU-sets},
put $f=|H|$ and choose the orthonormal normal frame
$e_3=\xi=H/f$, $e_4=\eta$. Assume that $\nabla f\neq0$,
$A_\xi$ is not a multiple of the identity, and $A_\eta\neq0$
on this neighborhood. Since the normal connection is flat,
\eqref{eq:Ricci-general} gives $[A_\xi,A_\eta]=0$. Thus a common
orthonormal principal frame $\{e_1,e_2\}$ gives
\begin{equation}\label{eq:rank-shape}
 \begin{aligned}
 h_{11}^{3}&=f(1+\lambda),& h_{22}^{3}&=f(1-\lambda),\\
 h_{11}^{4}&=f\mu,& h_{22}^{4}&=-f\mu,\\
 h_{12}^{3}&=h_{12}^{4}=0,&\lambda\mu&\neq0.
 \end{aligned}
\end{equation}
Here $\lambda$ and $\mu$ are smooth real-valued functions, and
$h_{ij}^{\alpha}=\langle h(e_i,e_j),e_\alpha\rangle$ for
$i,j=1,2$ and $\alpha=3,4$. Define
\begin{equation}\label{eq:rank-variables}
 u=\frac{e_1(f)}{f},\qquad v=\frac{e_2(f)}{f},\qquad
 r=\lambda^2+\mu^2.
\end{equation}
Since $\nabla f=f(ue_1+ve_2)$ and the eigenvalues in
\eqref{eq:rank-shape} are distinct, $\nabla f$ is a common principal
direction if and only if $uv=0$.

The tangent and normal connection coefficients are given by
\begin{equation}\label{eq:rank-frame}
 \omega_1^2\left(e_i\right)=\langle\nabla_{e_i}e_1,e_2\rangle,
 \qquad
 \omega_3^4\left(e_i\right)=\langle\nabla^\perp_{e_i}e_3,e_4\rangle,
 \quad i=1,2.
\end{equation}
The Lie bracket and Gauss equations become
\begin{align}
 [e_1,e_2]&=-\omega_1^2\left(e_1\right)e_1
             -\omega_1^2\left(e_2\right)e_2,
 \label{eq:rank-bracket}\\
 K&=e_2\left(\omega_1^2\left(e_1\right)\right)
       -e_1\left(\omega_1^2\left(e_2\right)\right)
       -\left(\omega_1^2\left(e_1\right)\right)^2
       -\left(\omega_1^2\left(e_2\right)\right)^2\notag\\
  &=c+f^2(1-r).
 \label{eq:rank-curvature}
\end{align}
Since the normal connection is flat, the connection form $\omega_3^4$
is closed. After restricting to a simply connected neighborhood, we can
write $\omega_3^4=\mathrm d\theta$ for a smooth function $\theta$.
Consequently,
\begin{align}\label{Fu1}
\nabla^\perp_{e_i}e_3=e_i(\theta)e_4,
\qquad
\nabla^\perp_{e_i}e_4=-e_i(\theta)e_3.
\end{align}
For later reference, we record the biharmonic equations for
$M^2$ in $N^4(c)$ in the frame \eqref{eq:rank-shape}. With the
Laplacian convention in Section~\ref{sec:preliminaries}, the normal
and tangential equations \eqref{eq:biharmonic-normal-general} and
\eqref{eq:biharmonic-tangent-general} take the form
\begin{subequations}\label{eq:surface-biharmonic-system}
\begin{align}
 &\Delta f
 =f\left(|\mathrm d\theta|^2+2f^2(1+\lambda^2)-2c\right),
 \label{eq:rank-normal-f}\\
 &\Delta\theta+2\langle\nabla\log f,\nabla\theta\rangle
 =2f^2\lambda\mu,
 \label{eq:rank-normal-theta}\\
 &(2+\lambda)e_1(f)+f\mu e_1(\theta)=0,
 \label{eq:surface-biharmonic-tangent-1}\\
 &(2-\lambda)e_2(f)-f\mu e_2(\theta)=0.
 \label{eq:surface-biharmonic-tangent-2}
\end{align}
\end{subequations}
The first two equations are the $e_3$- and $e_4$-components of the
normal equation; the last two are the $e_1$- and $e_2$-components
of the tangential equation after division by $2f$.
Since $H=fe_3$, we obtain
\[
\nabla^\perp_{e_i}H
=e_i(f)e_3+fe_i(\theta)e_4.
\]
Using \eqref{Fu1} to differentiate once more and taking the trace gives
\begin{align*}
 \langle\Delta^\perp H,e_3\rangle
 &=\Delta f-f|\mathrm d\theta|^2,\\
 \langle\Delta^\perp H,e_4\rangle
 &=f\Delta\theta+2\langle\nabla f,\nabla\theta\rangle.
\end{align*}
Moreover, \eqref{eq:rank-shape} yields
\[
 \sum_{i,j=1}^2(h_{ij}^3)^2=2f^2(1+\lambda^2),
 \qquad
 \sum_{i,j=1}^2h_{ij}^3h_{ij}^4=2f^2\lambda\mu.
\]
Substituting these identities into \eqref{eq:biharmonic-normal-general}
gives \eqref{eq:rank-normal-f} and \eqref{eq:rank-normal-theta},
where the $e_4$-component has been divided by $f>0$.
It follows from \eqref{eq:rank-shape} that
\begin{align*}
	A_{\nabla^\perp_{e_1}H}(e_1)
	&=\left(f(1+\lambda)e_1(f)
	+f^2\mu e_1(\theta)\right)e_1,\\
	A_{\nabla^\perp_{e_2}H}(e_2)
	&=\left(f(1-\lambda)e_2(f)
	-f^2\mu e_2(\theta)\right)e_2.
\end{align*}
On the other hand,
\[
\operatorname{grad}|H|^2
=2f\left(e_1(f)e_1+e_2(f)e_2\right).
\]
Substituting the above expressions into the tangential
biharmonic equation \eqref{eq:biharmonic-tangent-general}, we obtain
\begin{align*}
	\Big((2+\lambda)e_1(f)
	+f\mu e_1(\theta)\Big)e_1+
	\Big((2-\lambda)e_2(f)
	-f\mu e_2(\theta)\Big)e_2=0.
\end{align*}
Since $\{e_1,e_2\}$ is an orthonormal frame, the two components
vanish separately, giving \eqref{eq:surface-biharmonic-tangent-1}
and \eqref{eq:surface-biharmonic-tangent-2}. Using
\eqref{eq:rank-variables} and $\mu\neq0$, we therefore obtain

\begin{equation}\label{eq:rank-alpha}
 e_1(\theta)=-\frac{2+\lambda}{\mu}u,
 \qquad e_2(\theta)=\frac{2-\lambda}{\mu}v.
\end{equation}
Substituting \eqref{eq:rank-shape} and \eqref{eq:rank-alpha}
into \eqref{eq:Codazzi-general}, we obtain
\begin{subequations}\label{eq:rank-Codazzi}
\begin{align}
 e_1(\lambda)&=-(1+2\lambda)u-2\lambda\omega_1^2\left(e_2\right),
 \label{eq:rank-Codazzi-1}\\
 e_2(\lambda)&=(1-2\lambda)v+2\lambda\omega_1^2\left(e_1\right),
 \label{eq:rank-Codazzi-2}\\
 e_1(\mu)&=-\frac{2-\lambda-\lambda^2+\mu^2}{\mu}u
                  -2\mu\omega_1^2\left(e_2\right),
 \label{eq:rank-Codazzi-3}\\
 e_2(\mu)&=-\frac{2+\lambda-\lambda^2+\mu^2}{\mu}v
                  +2\mu\omega_1^2\left(e_1\right).
 \label{eq:rank-Codazzi-4}
\end{align}
\end{subequations}

\begin{lemma}\label{lem:flat-first-integral}
	On an open set where $uv\neq0$, the function
	\begin{equation}\label{eq:rank-m}
		\beta=\frac{\lambda(4-r)}{\mu r}
	\end{equation}
	depends locally only on $f$.
\end{lemma}

\begin{proof}
	Set
	\begin{equation}\label{eq:rank-XYz}
		X=\frac{e_1}{u},\qquad
		Y=\frac{e_2}{v},\qquad
		z=\frac{\mu}{\lambda},
	\end{equation}
	and write
	\begin{equation}\label{eq:rank-PQ-definition}
		P=\frac{\omega_1^2\left(e_1\right)}{v},\qquad
		Q=\frac{\omega_1^2\left(e_2\right)}{u}.
	\end{equation}
	The Codazzi equations
	\eqref{eq:rank-Codazzi-1}--\eqref{eq:rank-Codazzi-4}
	become
	\begin{equation}\label{eq:rank-directional-coefficients}
		\begin{aligned}
			X\left(\lambda\right)
			&=-1-2\lambda-2\lambda Q,
			&Y\left(\lambda\right)
			&=1-2\lambda+2\lambda P,\\
			X\left(\mu\right)
			&=-\frac{2-\lambda-\lambda^2+\mu^2}{\mu}-2\mu Q,
			&Y\left(\mu\right)
			&=-\frac{2+\lambda-\lambda^2+\mu^2}{\mu}+2\mu P.
		\end{aligned}
	\end{equation}
	Differentiating $z=\mu/\lambda$ and using
	\eqref{eq:rank-directional-coefficients}, we obtain
	\begin{equation}\label{eq:rank-k}
		X\left(z\right)
		=\frac{(\lambda+1)r-2\lambda}{\mu\lambda^2},
		\qquad
		Y\left(z\right)
		=\frac{(\lambda-1)r-2\lambda}{\mu\lambda^2}.
	\end{equation}
	By \eqref{eq:rank-variables}, we have
	$X\left(f\right)=Y\left(f\right)=f$. Hence
	\[
	\left(X-Y\right)\left(f\right)=0,
	\qquad
	\left(X-Y\right)\left(z\right)
	=\frac{2r}{\mu\lambda^2}\neq0.
	\]
	Combining \eqref{eq:rank-alpha} with \eqref{eq:rank-k}
	and using \eqref{eq:rank-m}, we obtain
	\[
	X\left(\theta+2\arctan z\right)
	=Y\left(\theta+2\arctan z\right)
	=-\beta.
	\]
	Since $X-Y$ spans $\ker\mathrm df$, the function
	$\theta+2\arctan z$ depends locally only on $f$.
	The identity $X\left(f\right)=f$ then shows that $\beta$
	also depends locally only on $f$.
\end{proof}

For the subsequent calculations, retain the notation above and put
\begin{equation}\label{eq:rank-gamma}
	\gamma=f\frac{\mathrm d\beta}{\mathrm df}.
\end{equation}
Since $\beta z=(4-r)/r$ and $r=\lambda^2(1+z^2)$, we have
\begin{equation}\label{eq:rank-r-s}
	r=\frac4{1+\beta z},\qquad
	\lambda^2=\frac{r}{1+z^2}.
\end{equation}
Differentiating $r=\lambda^2+\mu^2$ and using
\eqref{eq:rank-directional-coefficients} yields
\begin{equation}\label{eq:rank-r-derivatives}
	X\left(r\right)=-2(r+2)-4rQ,\qquad
	Y\left(r\right)=-2(r+2)+4rP.
\end{equation}
Moreover, \eqref{eq:rank-gamma} gives
$X\left(\beta\right)=Y\left(\beta\right)=\gamma$.
Differentiating $\beta z=4/r-1$ along $X$ and $Y$,
and substituting \eqref{eq:rank-r-derivatives}, we obtain
\begin{equation}\label{eq:rank-beta-derivatives}
	\begin{aligned}
		z\gamma+\beta X\left(z\right)
		&=\frac{8(r+2)+16rQ}{r^2},\\
		z\gamma+\beta Y\left(z\right)
		&=\frac{8(r+2)-16rP}{r^2}.
	\end{aligned}
\end{equation}
Subtracting the second equation from the first, and then adding
the two equations, we obtain, respectively, by
\eqref{eq:rank-k} and \eqref{eq:rank-m},
\begin{equation}\label{eq:rank-PQ}
	\begin{aligned}
		P+Q
		&=\frac{r(4-r)}{8\lambda\mu^2},\\
		P-Q
		&=\frac{r+2}{r}
		-\frac{(4-r)(r-2)}{8\mu^2}
		-\frac{r\mu}{8\lambda}\gamma.
	\end{aligned}
\end{equation}
Equation~\eqref{eq:rank-PQ} expresses $P$ and $Q$ as rational functions of
$\lambda,\mu,\gamma$, where $r=\lambda^2+\mu^2$.
A prime on $\gamma$ denotes differentiation with respect to $f$.
By the definitions of $X$ and $Y$, we have
\[
X(f)=Y(f)=f.
\]
Consequently,
\[
X(\gamma)=Y(\gamma)=f\gamma'.
\]
Partial derivatives such as $P_\lambda$ and $Q_\mu$
are taken with the other variables fixed. Since $\gamma$
depends only on $f$, the chain rule gives
\begin{equation}\label{eq:rank-chain-rule}
 X\left(P\right)=P_{\lambda}X\left(\lambda\right)+P_\mu X\left(\mu\right)+P_\gamma f\gamma',
 \qquad
 Y\left(P\right)=P_{\lambda}Y\left(\lambda\right)+P_\mu Y\left(\mu\right)+P_\gamma f\gamma'.
\end{equation}
The same chain rule applies to $Q$.

\begin{lemma}\label{lem:flat-compatibility}
On the open set where $uv\neq0$, the mixed derivatives of $u,v$
are given by
\begin{subequations}\label{eq:rank-cross}
\begin{align}
 e_2(u)&=\left(2P-Q-\frac{3\mu^2+4-\lambda^2}{2\mu^2}\right)uv,
 \label{eq:rank-cross-u}\\
 e_1(v)&=\left(P-2Q-\frac{3\mu^2+4-\lambda^2}{2\mu^2}\right)uv.
 \label{eq:rank-cross-v}
\end{align}
\end{subequations}
The remaining first derivatives are given by
\begin{equation}\label{eq:rank-normal-u}
 \begin{aligned}
 e_1(u)={}&\left(\frac{(\lambda+2)^2}{4\mu^2}-\frac74-2Q\right)u^2
       +\left(\frac{3(\lambda-2)^2}{4\mu^2}+\frac34-P\right)v^2\\
       &+\left(1+\lambda^2-\frac{\lambda(1+r)}2\right)f^2-\frac{2-\lambda}{2}c,
 \end{aligned}
\end{equation}
\begin{equation}\label{eq:rank-normal-v}
 \begin{aligned}
 e_2(v)={}&\left(\frac{3(\lambda+2)^2}{4\mu^2}+\frac34+Q\right)u^2
       +\left(\frac{(\lambda-2)^2}{4\mu^2}-\frac74+2P\right)v^2\\
       &+\left(1+\lambda^2+\frac{\lambda(1+r)}2\right)f^2-\frac{2+\lambda}{2}c.
 \end{aligned}
\end{equation}
Moreover, the Gauss equation determines the remaining derivative
$f\gamma'$ by
\begin{equation}\label{eq:rank-n1}
 \begin{aligned}
 \frac{\mu r(u^2+v^2)}{16\lambda}f\gamma'
 ={}&\left(P_{\lambda}Y\left(\lambda\right)+P_\mu Y\left(\mu\right)-P^2\right)v^2\\
    &-\left(Q_{\lambda}X\left(\lambda\right)+Q_\mu X\left(\mu\right)+Q^2\right)u^2\\
    &+P e_2(v)-Q e_1(u)-c-f^2(1-r).
 \end{aligned}
\end{equation}
Every quantity on the right-hand side is given by
\eqref{eq:rank-PQ}, \eqref{eq:rank-directional-coefficients},
\eqref{eq:rank-normal-u}, and \eqref{eq:rank-normal-v}.
\end{lemma}

\begin{proof}
We use the  biharmonic system
\eqref{eq:surface-biharmonic-system} established above and verify
the claimed derivative identities in stages.
Applying \eqref{eq:rank-bracket} to $\log f$ and using
\eqref{eq:rank-variables} and \eqref{eq:rank-PQ-definition}, we obtain
\begin{equation}\label{eq:rank-logf-commutator}
 e_1(v)-e_2(u)=-(P+Q)uv.
\end{equation}
Applying \eqref{eq:rank-bracket} to $\theta$ and substituting
\eqref{eq:rank-alpha} and \eqref{eq:rank-Codazzi}, we obtain
\begin{equation}\label{eq:rank-theta-commutator}
 \begin{aligned}
 &(2+\lambda)e_2(u)+(2-\lambda)e_1(v)\\
 &\quad=\left((6+\lambda)P+(\lambda-6)Q
       -\frac{2(3\mu^2+4-\lambda^2)}{\mu^2}\right)uv.
 \end{aligned}
\end{equation}
For clarity, put $U=e_2(u)/(uv)$ and $V=e_1(v)/(uv)$.
The two equations form the $2\times2$ system
\[
 \begin{aligned}
 V-U&=-(P+Q),\\
 (2+\lambda)U+(2-\lambda)V
 &=(6+\lambda)P+(\lambda-6)Q
 -\frac{2(3\mu^2+4-\lambda^2)}{\mu^2},
 \end{aligned}
\]
whose determinant is $-4$, independently of $\lambda$.
Solving for $U,V$ and multiplying by $uv$ yields \eqref{eq:rank-cross-u} and
\eqref{eq:rank-cross-v}.

We next use the two normal equations
\eqref{eq:rank-normal-f} and \eqref{eq:rank-normal-theta}.
Equations \eqref{eq:rank-variables} and
\eqref{eq:rank-PQ-definition} give
\begin{equation}\label{eq:rank-laplacian-f}
 \frac{\Delta f}{f}
 =e_1(u)+e_2(v)+(1+Q)u^2+(1-P)v^2,
\end{equation}
and substituting \eqref{eq:rank-alpha} and
\eqref{eq:rank-laplacian-f} into \eqref{eq:rank-normal-f} gives
\begin{equation}\label{eq:rank-normal-sum}
 \begin{aligned}
 e_1(u)+e_2(v)={}&
 \left(\frac{(\lambda+2)^2}{\mu^2}-1-Q\right)u^2
 +\left(\frac{(\lambda-2)^2}{\mu^2}-1+P\right)v^2\\
 &+2f^2(1+\lambda^2)-2c,
 \end{aligned}
\end{equation}
whereas substitution of \eqref{eq:rank-alpha},
\eqref{eq:rank-Codazzi}, and \eqref{eq:rank-PQ-definition} into
\eqref{eq:rank-normal-theta} yields
\begin{align}
 -(2+\lambda)e_1(u)+(2-\lambda)e_2(v)
 ={}&\left((\lambda+6)Q+\lambda+5
       -\frac{(\lambda-1)(\lambda+2)^2}{\mu^2}\right)u^2\label{eq:rank-normal-difference}\\
 &+\left((6-\lambda)P+\lambda-5
       -\frac{(\lambda+1)(\lambda-2)^2}{\mu^2}\right)v^2
       +2f^2\lambda\mu^2.\nonumber
 \end{align}
To check the remaining algebra, write
$A=e_1(u)$, $B=e_2(v)$ and denote the right-hand sides
of \eqref{eq:rank-normal-sum} and
\eqref{eq:rank-normal-difference} by $S$ and $T$.
The coefficient matrix is
\[
 \begin{pmatrix}1&1\\-(2+\lambda)&2-\lambda\end{pmatrix},
 \qquad \det=4,
\]
so that $A=((2-\lambda)S-T)/4$ and
$B=((2+\lambda)S+T)/4$.
In particular, the $c$-terms are $-(2-\lambda)c/2$
and $-(2+\lambda)c/2$, respectively; the $f^2$-terms
simplify using $r=\lambda^2+\mu^2$.
Multiply \eqref{eq:rank-normal-sum} by $2-\lambda$ and subtract
\eqref{eq:rank-normal-difference}. Division by $4$ gives
\eqref{eq:rank-normal-u}. Multiplying \eqref{eq:rank-normal-sum} by
$2+\lambda$ and adding \eqref{eq:rank-normal-difference} gives
\eqref{eq:rank-normal-v} after the same division.

Substituting $\omega_1^2\left(e_1\right)=Pv$ and
$\omega_1^2\left(e_2\right)=Qu$ from \eqref{eq:rank-PQ-definition}
into \eqref{eq:rank-curvature}, we obtain
\begin{equation}\label{eq:rank-Gauss}
 (Y\left(P\right)-P^2)v^2-(X\left(Q\right)+Q^2)u^2+P e_2(v)-Q e_1(u)-c-f^2(1-r)=0.
\end{equation}
By \eqref{eq:rank-PQ}, we get
\[
 P_\gamma=-\frac{r\mu}{16\lambda},\qquad
 Q_\gamma=\frac{r\mu}{16\lambda}.
\]
Consequently the terms involving $f\gamma'$ in
\eqref{eq:rank-Gauss} are exactly
\[
 \left(P_\gamma v^2-Q_\gamma u^2\right)f\gamma'
 =-\frac{\mu r}{16\lambda}(u^2+v^2)f\gamma'.
\]
Using the chain rule \eqref{eq:rank-chain-rule} in
\eqref{eq:rank-Gauss} and moving the term containing $f\gamma'$ to the
other side gives \eqref{eq:rank-n1}. The coefficient of $f\gamma'$ is
nonzero because $\lambda\mu\neq0$, $r>0$, and $u^2+v^2>0$.
\end{proof}

\begin{lemma}\label{lem:flat-normal-principal}
Let $M^2$ be a biharmonic surface with flat normal bundle in $N^4(c)$,
where $c\leq0$. On an open set on which $H=f\xi$, $f>0$, and
$\nabla f\neq0$, suppose that $A_\xi$ is not a multiple of the identity
and that $A_\eta\neq0$, where $\eta$ is a unit normal orthogonal to
$\xi$. Then $\nabla f$ is a common principal direction of $A_\xi$ and
$A_\eta$.
\end{lemma}

\begin{proof}
Suppose that the conclusion fails. Restricting to a nonempty connected
open set, we may assume $uv\neq0$. We use Lemmas~\ref{lem:flat-first-integral}
and \ref{lem:flat-compatibility}. Define
\begin{equation}\label{eq:rank-GZ}
 \mathcal G=\frac{u^2+v^2}{f^2},\qquad
 \mathcal Z=\frac{\lambda(u^2-v^2)}{f^2}.
\end{equation}
Then $\mathcal G=|\nabla f|^2/f^4$, and
\begin{equation}\label{eq:rank-uv-GZ}
 \frac{u^2}{f^2}=\frac12\left(\mathcal G+\frac{\mathcal Z}{\lambda}\right),
 \qquad
 \frac{v^2}{f^2}=\frac12\left(\mathcal G-\frac{\mathcal Z}{\lambda}\right).
\end{equation}

Applying \eqref{eq:rank-bracket} to $u$ and $v$ gives, respectively,
\begin{subequations}
\begin{align}
 e_2\left(e_1(u)\right)-e_1\left(e_2(u)\right)
 -Pv\,e_1(u)-Qu\,e_2(u)&=0,
 \label{eq:rank-commutator-u}\\
 e_2\left(e_1(v)\right)-e_1\left(e_2(v)\right)
 -Pv\,e_1(v)-Qu\,e_2(v)&=0.
 \label{eq:rank-commutator-v}
\end{align}
\end{subequations}
In \eqref{eq:rank-commutator-u} and
\eqref{eq:rank-commutator-v}, substitute the first derivatives from
\eqref{eq:rank-cross-u}, \eqref{eq:rank-cross-v},
\eqref{eq:rank-normal-u}, and \eqref{eq:rank-normal-v}.
Differentiating the coefficients by \eqref{eq:rank-chain-rule} and
\eqref{eq:rank-directional-coefficients}, we eliminate $f\gamma'$
using \eqref{eq:rank-n1}. Divide the resulting form of
\eqref{eq:rank-commutator-u} by $vf^2$ and the resulting form of
\eqref{eq:rank-commutator-v} by $uf^2$.

Adding these normalized equations and multiplying by
$-2\lambda z^3(1+\beta z)$ gives \eqref{eq:rank-linear-first} after
substitution of \eqref{eq:rank-uv-GZ} and \eqref{eq:rank-r-s}.
Subtracting the normalized form of \eqref{eq:rank-commutator-v}
from that of \eqref{eq:rank-commutator-u}, multiplying by
$4z^3(1+\beta z)/(1+z^2)$, and making the same substitutions gives
\eqref{eq:rank-linear-second}. In both calculations,
$\mu=\lambda z$ from \eqref{eq:rank-XYz}, so that
\eqref{eq:rank-r-s} eliminates $\lambda^2$ and $\mu$.
Collecting the coefficients of $\mathcal G$ and $\mathcal Z$, we obtain
\begin{equation}\label{eq:rank-linear-first}
\begin{aligned}
 &\bigl[\beta^2z(3z^4+10z^2+11)
       +2\beta(z^4+7)+4\gamma z^2-12z\bigr]\mathcal G\\
 &\quad-(1+z^2)\bigl[\beta^2z(z^2-9)
       +\beta(9z^2-7)+2\gamma z^2+6z\bigr]\mathcal Z\\
 &\quad =4\beta z^2(z^2+5)\frac{c}{f^2}-\frac{4z^2}{(1+z^2)(1+\beta z)^2}
       \Bigl[\beta^3z^2(1+z^2)(z^2+5)
       \\
       &\qquad\,+2\beta^2z(z^4+2z^2+15)+\beta(z^4-22z^2+17)-8\gamma z^2-12z\Bigr],
\end{aligned}
\end{equation}
\begin{equation}\label{eq:rank-linear-second}
\begin{aligned}
 &\bigl[3\beta^3z^2(z^2+3)+\beta^2z(13z^2+21)
       +6\beta\gamma z^3+2\beta(6z^2+7)+8\gamma z^2\bigr]\mathcal G\\
 &\quad +7\beta(1+z^2)(1+\beta z)^2\mathcal Z\\
 &\quad =\frac{4z^2}{1+z^2}
       \bigl[\beta^2z(z^2+3)+\beta(8z^2+5)
       +2\gamma z^2+5z\bigr]\frac{c}{f^2}\\
 &\qquad -\frac{4z^2}{(1+z^2)(1+\beta z)}
       \Bigl[\beta^3z^2(z^2+3)+3\beta^2z(3z^2+8)\\
 &\qquad \hphantom{{}-\frac{4z^2}{(1+z^2)(1+\beta z)}\Bigl[}
       {}+\beta(2\gamma z^3+13z^2+17)+2\gamma z^2+9z\Bigr].
\end{aligned}
\end{equation}
To check the elimination, one may treat
$(f,z)$ as local coordinates: indeed $df\ne0$ and
$(X-Y)f=0$, whereas $(X-Y)z=2r/(\mu\lambda^2)\ne0$.
In the substitutions above, $\beta,\gamma$ and $c/f^2$
are held fixed when collecting powers of $z$; no derivative
of these functions with respect to $z$ is taken.
By \eqref{eq:rank-XYz} and \eqref{eq:rank-r-s}, we have
$z=\mu/\lambda\neq0$ and $1+\beta z=4/r>0$. Hence the denominators
in \eqref{eq:rank-linear-first} and \eqref{eq:rank-linear-second}
are nonzero.

If $\beta=0$ on an open set, then \eqref{eq:rank-gamma} gives
$\gamma=0$. Equation
\eqref{eq:rank-linear-second} reduces to
\[
 0=\frac{4z^3}{1+z^2}\left(\frac{5c}{f^2}-9\right),
\]
which is impossible because $z\neq0$ and $c\leq0$.
We may therefore restrict to a smaller nonempty open set on which
$\beta\neq0$.

By Lemma~\ref{lem:flat-first-integral} and \eqref{eq:rank-gamma},
$\beta$ and $\gamma$ depend only on $f$. Fixing $f$ in the coordinates
$(f,z)$ therefore fixes the values of $\beta$, $\gamma$, and $c/f^2$.
The coefficients of $\mathcal G$ in
\eqref{eq:rank-linear-first} and \eqref{eq:rank-linear-second} have
leading terms $3\beta^2z^5$ and $3\beta^3z^4$, respectively.
Since $\beta\neq0$, neither coefficient is the zero polynomial.
Dividing by the respective coefficients of $\mathcal G$ and expanding
the resulting rational functions at infinity, we obtain
\begin{equation}\label{eq:rank-Laurent-coefficients}
\begin{aligned}
 \mathcal G+\left(-\frac13+O(z^{-1})\right)\mathcal Z
 &=\frac{4(c/f^2-1)}{3\beta z}
       -\frac{8(c/f^2-1)}{9\beta^2z^2}+O(z^{-3}),\\
 \mathcal G+\left(\frac73+O(z^{-1})\right)\mathcal Z
 &=\frac{4(c/f^2-1)}{3\beta z}
       +\frac{44(c/f^2-1)}{9\beta^2z^2}+O(z^{-3}).
\end{aligned}
\end{equation}
By \eqref{eq:rank-Laurent-coefficients}, the difference between
the coefficients of $\mathcal Z$ is $8/3+O(z^{-1})$.
Thus the determinant of the normalized equations is not the zero
rational function. For fixed $f$, the zeros of this determinant and
of the two original coefficients of $\mathcal G$ are isolated.
After restricting to a smaller coordinate rectangle on which these
quantities are nonzero, we obtain from Cramer's rule that
$\mathcal G$ and $\mathcal Z$ are rational functions of $z$.

Subtracting the first equation in \eqref{eq:rank-Laurent-coefficients}
from the second and substituting the result into the first, we obtain
\begin{align}
 \mathcal Z&=\frac{13(c/f^2-1)}{6\beta^2z^2}+O(z^{-3}),
 \label{eq:rank-z-expansion}\\
 \mathcal G&=\frac{4(c/f^2-1)}{3\beta z}
            -\frac{c/f^2-1}{6\beta^2z^2}+O(z^{-3}).
 \label{eq:rank-gz-expansion}
\end{align}

Since $X-Y$ is tangent to the level curves of $f$,
\eqref{eq:rank-k} and \eqref{eq:rank-r-s} give
\begin{equation}\label{eq:rank-level-derivative}
 \lambda(X-Y)=\frac{2(1+z^2)}{z}\frac{\partial}{\partial z},
 \qquad\text{with }f\text{ fixed}.
\end{equation}
Applying $\lambda(X-Y)$ to \eqref{eq:rank-GZ} and using
\eqref{eq:rank-XYz}, we obtain
\[
 \lambda\left(X-Y\right)\left(\mathcal G\right)
 =\frac{2\lambda}{f^2}
 \left(e_1(u)-e_2(v)+\frac vu e_1(v)-\frac uv e_2(u)\right).
\]
Substituting \eqref{eq:rank-cross-u}, \eqref{eq:rank-cross-v},
\eqref{eq:rank-normal-u}, and \eqref{eq:rank-normal-v} into this
equation, and then using \eqref{eq:rank-PQ},
\eqref{eq:rank-uv-GZ}, and \eqref{eq:rank-level-derivative}, we obtain
\begin{equation}\label{eq:rank-Lg}
\begin{aligned}
 \frac{2(1+z^2)}{z}\frac{\partial\mathcal G}{\partial z}
 ={}&\left(\frac{r-4}{2}+\frac{r-12}{2z^2}\right)\mathcal G
      -\frac{2(1+z^2)}{z^2}\mathcal Z\\
 &+\frac{2r}{1+z^2}\left(\frac{c}{f^2}-r-1\right),
 \qquad r=\frac4{1+\beta z}.
\end{aligned}
\end{equation}
For fixed $f$, both sides of \eqref{eq:rank-Lg} are rational
functions of $z$. Multiplying the difference between the two sides
by a common denominator gives a polynomial that vanishes on an open
interval. Hence this polynomial is identically zero, and the Laurent
coefficients of the two sides at infinity agree.

Differentiating \eqref{eq:rank-gz-expansion} and substituting into
the left-hand side of \eqref{eq:rank-Lg}, we obtain
\[
 -\frac{8(c/f^2-1)}{3\beta z}
 +\frac{2(c/f^2-1)}{3\beta^2z^2}+O(z^{-3}).
\]
Substituting \eqref{eq:rank-z-expansion},
\eqref{eq:rank-gz-expansion}, and \eqref{eq:rank-r-s} into the
right-hand side of \eqref{eq:rank-Lg} gives
\[
 -\frac{8(c/f^2-1)}{3\beta z}
 -\frac{4(c/f^2-1)}{3\beta^2z^2}+O(z^{-3}).
\]
Comparison of the coefficients of $z^{-2}$ in \eqref{eq:rank-Lg}
therefore gives
\[
 \frac{2(c/f^2-1)}{\beta^2}=0.
\]
Since $c\leq0$ and $f>0$, we have $c/f^2-1<0$, which is a
contradiction. Thus $uv=0$, and $\nabla f$ is a common principal
direction of the two shape operators.
\end{proof}

\section{Minimality of surfaces with flat normal bundle}\label{sec:vector}

\begin{proof}[Reduction of Theorem~\ref{thm:flat} to
Proposition~\ref{lem:vector-ode-rigidity}]
Suppose that $\mathcal V=\{x\in M:H(x)\neq0\}$ is nonempty.
On a sufficiently small open subset of $\mathcal V$, choose
$e_3=\xi=H/f$, $e_4=\eta$, where $f=|H|>0$.
The $e_3$-component of \eqref{eq:biharmonic-normal-general} gives
\begin{equation}\label{eq:flat-final-radial}
 \Delta f=f\left(
 \sum_{i=1}^2\left(\omega_3^4\left(e_i\right)\right)^2
 +\sum_{i,j=1}^2\left(h_{ij}^{3}\right)^2-2c\right).
\end{equation}
Since $h_{11}^{3}+h_{22}^{3}=2f$, equation~\eqref{eq:flat-final-radial}
and the Cauchy--Schwarz inequality yield $\Delta f\geq2f(f^2-c)>0$.
It follows that $f$ is not constant on any nonempty open subset of
$\mathcal V$. Restricting to a smaller nonempty connected open set,
we may assume that $\nabla f\neq0$.

Suppose first that $A_\xi=f\mathrm{Id}$ on an open set. At each point, choose a
basis diagonalizing $A_\eta$, and write the eigenvalues as
$f\mu,-f\mu$. Equations \eqref{eq:Codazzi-general} and
\eqref{eq:biharmonic-tangent-general} give, respectively,
\[
 \begin{aligned}
 &e_1(f)=-f\mu\omega_3^4\left(e_1\right),&
 &e_2(f)= f\mu\omega_3^4\left(e_2\right),\\
 &2e_1(f)+f\mu\omega_3^4\left(e_1\right)=0,&
 &2e_2(f)-f\mu\omega_3^4\left(e_2\right)=0.
 \end{aligned}
\]
Substituting the first row into the second yields $e_1(f)=e_2(f)=0$,
contrary to $\nabla f\neq0$. After restriction to a smaller nonempty
open set, $A_\xi$ has two distinct eigenvalues.

Next suppose that $A_\eta=0$ on this open set. Write the eigenvalues
of $A_\xi$ as $f(1+\lambda),f(1-\lambda)$, where $\lambda\neq0$. The tangential
biharmonic equation \eqref{eq:biharmonic-tangent-general} and the
$\eta$-component of \eqref{eq:Codazzi-general} give
\[
 \begin{aligned}
 (2+\lambda)e_1(f)&=0,& (2-\lambda)e_2(f)&=0,\\
 (1-\lambda)\omega_3^4\left(e_1\right)&=0,&
 (1+\lambda)\omega_3^4\left(e_2\right)&=0.
 \end{aligned}
\]
Since $\nabla f\neq0$, the first row gives $\lambda=2$ or $\lambda=-2$.
The second row then yields
$\omega_3^4\left(e_1\right)=\omega_3^4\left(e_2\right)=0$.
Hence $\nabla^\perp\xi=0$. Since $f>0$, this contradicts
Theorem~\ref{thm:main}.
We may therefore restrict to a nonempty open set on which $\lambda\mu\neq0$.

Using Lemma~\ref{lem:flat-normal-principal}, we may choose the principal
frame so that
\[
 e_1=\frac{\nabla f}{|\nabla f|},\qquad e_2(f)=0,
 \qquad u=\frac{e_1(f)}{f}>0.
\]
The second fundamental form is given by \eqref{eq:rank-shape}.
Substituting $e_2(f)=0$ into \eqref{eq:rank-alpha},
\eqref{eq:rank-Codazzi-2}, and \eqref{eq:rank-Codazzi-4}, and applying
\eqref{eq:rank-bracket} to $f$, we obtain
\begin{equation}\label{eq:flat-final-transverse}
 \begin{aligned}
 \omega_3^4\left(e_1\right)&=-\frac{2+\lambda}{\mu}u,&
 \omega_3^4\left(e_2\right)&=0,\\
 e_2(\lambda)&=2\lambda\omega_1^2\left(e_1\right),&
 e_2(\mu)&=2\mu\omega_1^2\left(e_1\right),\\
 e_2(u)&=u\omega_1^2\left(e_1\right).
 \end{aligned}
\end{equation}
Since $\mathrm d\omega_3^4=0$, equations
\eqref{eq:rank-bracket} and \eqref{eq:flat-final-transverse} give
\begin{equation}\label{eq:flat-normal-closed}
 e_2\left(\omega_3^4\left(e_1\right)\right)
 =\omega_1^2\left(e_1\right)\omega_3^4\left(e_1\right)
 =-\frac{2+\lambda}{\mu}u\omega_1^2\left(e_1\right).
\end{equation}
Differentiating the first equation in \eqref{eq:flat-final-transverse}
in the direction $e_2$ and substituting the remaining equations gives
\[
 \begin{aligned}
 e_2\left(-\frac{2+\lambda}{\mu}u\right)
 &=-\frac{e_2(\lambda)}{\mu}u
   -\frac{2+\lambda}{\mu}e_2(u)
   +\frac{2+\lambda}{\mu^2}u e_2(\mu)\\
 &=\frac{2-\lambda}{\mu}u\omega_1^2\left(e_1\right).
 \end{aligned}
\]
Comparison with \eqref{eq:flat-normal-closed} yields
$4u\omega_1^2\left(e_1\right)/\mu=0$. Since $u>0$ and $\mu\neq0$,
equation~\eqref{eq:flat-final-transverse} therefore gives
\begin{equation}\label{eq:flat-final-one-variable}
 \omega_1^2\left(e_1\right)=0,
 \qquad e_2(\lambda)=e_2(\mu)=e_2(u)=0.
\end{equation}
By \eqref{eq:flat-final-one-variable}, the vector field $e_1$ is
geodesic, and $f,\lambda,\mu,u$ are constant along the local level curves
of $f$. Equation~\eqref{eq:rank-Codazzi-1} gives
\[
 \omega_1^2\left(e_2\right)
 =-\frac{e_1(\lambda)+(1+2\lambda)u}{2\lambda}.
\]
Applying $e_2$ and using \eqref{eq:rank-bracket} and
\eqref{eq:flat-final-one-variable}, we obtain
$e_2\left(\omega_1^2\left(e_2\right)\right)=0$.

Choose coordinates $(t,y)$ with $e_1=\partial_t$ and
$g=\mathrm dt^2+\chi(t,y)^2\mathrm dy^2$. Then
$\partial_t\log \chi=\omega_1^2\left(e_2\right)$ depends only on $t$, and hence
$\chi(t,y)=\rho(t)\chi_0(y)$ for positive functions $\rho,\chi_0$.
Replacing $y$ by $\int \chi_0(y)\,\mathrm dy$, we obtain
\begin{equation}\label{eq:flat-final-metric}
 g=\mathrm dt^2+\rho(t)^2\mathrm dy^2,
 \qquad b=-\omega_1^2\left(e_2\right),
 \qquad \rho'=-b\rho,\quad\rho>0.
\end{equation}
A prime denotes differentiation with respect to $t$.
Since $\omega_3^4\left(e_2\right)=0$, the rotation from $\{\xi,\eta\}$ to a
parallel normal frame depends only on $t$. Choose such a frame
$\{\nu_1,\nu_2\}$ on a smaller simply connected open set, and put
\begin{equation}\label{eq:flat-final-components}
 H_\alpha=\frac12\left(h_{11}^{\alpha}+h_{22}^{\alpha}\right),
 \qquad k_\alpha=h_{22}^{\alpha},\qquad \alpha=1,2.
\end{equation}
Since $\{\nu_1,\nu_2\}$ is an orthonormal normal frame, we have
\[
H=\sum_{\alpha=1}^{2}H_\alpha\nu_\alpha,
\qquad
f=|H|=\left(\sum_{\alpha=1}^{2}H_\alpha^2\right)^{1/2},
\]
where $f$ is the same mean curvature function introduced earlier and
these functions depend only on $t$. With this frame,
\begin{equation}\label{eq:flat-final-shape-parallel}
 h_{11}^{\alpha}=2H_\alpha-k_\alpha,
 \qquad h_{22}^{\alpha}=k_\alpha,
 \qquad h_{12}^{\alpha}=0.
\end{equation}

Substituting \eqref{eq:flat-final-shape-parallel} into
\eqref{eq:Codazzi-general} gives
$k_\alpha'=2b(k_\alpha-H_\alpha)$. The Gauss equation
\eqref{eq:Gauss-general}, together with $K=b'-b^2$, gives
\[
 b'=b^2+c+\sum_{\alpha=1}^2
             \left(2H_\alpha k_\alpha-k_\alpha^2\right).
\]
Thus \eqref{eq:vo-G} and \eqref{eq:vo-b} hold. Put
\[
 f^2=\sum_{\alpha=1}^2H_\alpha^2,
 \qquad d=\sum_{\alpha=1}^2H_\alpha(k_\alpha-H_\alpha).
\]
In the parallel normal frame, \eqref{eq:flat-final-metric} gives
$\langle\Delta^\perp H,\nu_\alpha\rangle=H_\alpha''-bH_\alpha'$.
Moreover, \eqref{eq:flat-final-shape-parallel} yields
\begin{align*}
 \sum_{i=1}^2\langle h(A_He_i,e_i),\nu_\alpha\rangle
 &=\sum_{\beta=1}^2H_\beta
       (2H_\beta-k_\beta)(2H_\alpha-k_\alpha)
       +\sum_{\beta=1}^2H_\beta k_\beta k_\alpha\\
 &=2f^2H_\alpha+2d(k_\alpha-H_\alpha).
\end{align*}
Thus \eqref{eq:biharmonic-normal-general} becomes \eqref{eq:vo-H}.
The $e_1$-component of \eqref{eq:biharmonic-tangent-general} gives
\[
 2\sum_{\alpha=1}^2(2H_\alpha-k_\alpha)H_\alpha'+(f^2)'=0.
\]
Since $(f^2)'=2\sum_{\alpha=1}^2H_\alpha H_\alpha'$, this reduces to
\eqref{eq:vo-tangent}. The $e_2$-component vanishes since all coefficients
depend only on $t$. Thus \eqref{eq:vo-G}--\eqref{eq:vo-tangent} hold with
$p=2$. Proposition~\ref{lem:vector-ode-rigidity} gives $H_1=H_2=0$,
contrary to $f>0$.
\end{proof}

In the following proposition, all sums over $\alpha$ range from $1$ to $p$.

\Needspace{19\baselineskip}
\begin{proposition}\label{lem:vector-ode-rigidity}
Let $c\leq0$ and $p\geq1$. Suppose that $H_\alpha,k_\alpha,b$ are
smooth real functions on a nonempty connected open interval $I$.
Define $d=\sum_\alpha H_\alpha(k_\alpha-H_\alpha)$.
If
\begin{align}
 k_\alpha'&=2b(k_\alpha-H_\alpha),
 \label{eq:vo-G}\\
 b'&=b^2+c+\sum_\alpha(2H_\alpha k_\alpha-k_\alpha^2),
 \label{eq:vo-b}\\
 H_\alpha''&=bH_\alpha'+2(f^2-c)H_\alpha
                         +2d(k_\alpha-H_\alpha),
 \label{eq:vo-H}\\
 0&=\sum_\alpha(3H_\alpha-k_\alpha)H_\alpha',
 \label{eq:vo-tangent}
\end{align}
then $H_\alpha\equiv0$ for every $\alpha$.
\end{proposition}

\begin{proof}
The function $f$ is smooth where $f>0$. Multiplying \eqref{eq:vo-H} by
$H_\alpha$, summing over $\alpha$, and using
$ff'=\sum_\alpha H_\alpha H_\alpha'$ and the Cauchy--Schwarz inequality,
we obtain
\begin{equation}\label{eq:vo-radial}
 f''-bf'
 =\frac{\sum_\alpha(H_\alpha')^2-(f')^2}{f}
        +2f(f^2-c)+\frac{2d^2}{f}
 \geq2f(f^2-c).
\end{equation}
Since $c\leq0$, equation~\eqref{eq:vo-radial} excludes a positive
constant value of $f$ on an open interval. Suppose that some $H_\alpha$
does not vanish identically. Restriction to a smaller interval,
translation of the parameter and, if necessary, reversal of $t$,
with $b(t)$ replaced by $-b(-t)$, give $f(0)>0$ and $f'(0)>0$.

Since $f^2=\sum_\alpha H_\alpha^2$, equations
\eqref{eq:vo-G}--\eqref{eq:vo-H} form a polynomial first-order system
for $k_\alpha,b,H_\alpha,H_\alpha'$. Let $t_*$ be the right endpoint
of the maximal interval of existence. The solution is real analytic,
so the identity theorem extends \eqref{eq:vo-tangent} from the initial
interval to the maximal interval.

Choose $\rho(0)>0$ so that
$\left(b(0)^2+\sum_\alpha k_\alpha(0)^2+c\right)\rho(0)^2\leq1$,
and define $\rho$ by
\begin{equation}\label{eq:vo-rho}
 \rho'=-b\rho,\qquad\rho(0)>0.
\end{equation}
Equation~\eqref{eq:vo-rho} gives $\rho>0$. Multiplying
\eqref{eq:vo-radial} by $\rho$ and using \eqref{eq:vo-rho}, we obtain
\begin{equation}\label{eq:vo-q-positive}
 (f'\rho)'\geq2f(f^2-c)\rho>0
 \qquad\text{as long as }f>0.
\end{equation}
It follows that $f$ and $f'$ remain positive, and
\begin{equation}\label{eq:vo-positive}
 f(t)\geq f(0)>0,\qquad
 f'(t)\rho(t)\geq\delta:=f'(0)\rho(0)>0,
 \qquad 0\leq t<t_*.
\end{equation}

Equations \eqref{eq:vo-G} and \eqref{eq:vo-b} give
\[
 \left(b^2+\sum_\alpha k_\alpha^2+c\right)'
 =2b\left(b^2+\sum_\alpha k_\alpha^2+c\right).
\]
Using \eqref{eq:vo-rho} and integrating, we obtain
\begin{equation}\label{eq:vo-first-integral}
 \left(b^2+\sum_\alpha k_\alpha^2+c\right)\rho^2=E,
\end{equation}
where $E$ is constant.
By the choice of $\rho(0)$, we have $E\leq1$. Put $c_0=\sqrt{-c}$.
Equation~\eqref{eq:vo-first-integral} yields
\begin{equation}\label{eq:vo-energy-bounds}
 |b|\rho\leq1+c_0\rho,
 \qquad
 \rho\left(\sum_\alpha k_\alpha^2\right)^{1/2}\leq1+c_0\rho.
\end{equation}
Equation \eqref{eq:vo-tangent} gives
\begin{equation}\label{eq:vo-tangent-product}
 \sum_\alpha(k_\alpha-H_\alpha)H_\alpha'
 =2\sum_\alpha H_\alpha H_\alpha'=2ff'.
\end{equation}
Differentiating $d$ and substituting \eqref{eq:vo-G} and
\eqref{eq:vo-tangent-product}, we obtain
\[
 \begin{aligned}
 d'&=\sum_\alpha H_\alpha'(k_\alpha-H_\alpha)
       +\sum_\alpha H_\alpha(k_\alpha'-H_\alpha')\\
   &=2ff'+2bd-ff'=2bd+ff'.
 \end{aligned}
\]
It follows from \eqref{eq:vo-rho} and \eqref{eq:vo-positive} that
\begin{equation}\label{eq:vo-Uderivative}
 (d\rho^2)'=ff'\rho^2>0.
\end{equation}
Next multiply \eqref{eq:vo-H} by $2H_\alpha'$ and sum over $\alpha$.
Using \eqref{eq:vo-rho} and \eqref{eq:vo-tangent-product}, we obtain
\[
 \left(\rho^2\sum_\alpha(H_\alpha')^2\right)'
 =4ff'(f^2-c)\rho^2+8ff'd\rho^2.
\]
By \eqref{eq:vo-Uderivative}, we also have
\[
 \left(4(f^2-c)d\rho^2\right)'
 =8ff'd\rho^2+4(f^2-c)ff'\rho^2.
\]
Subtracting these two equations and integrating, we obtain
\begin{equation}\label{eq:vo-second-integral}
 \rho^2\sum_\alpha(H_\alpha')^2
       -4(f^2-c)d\rho^2=C
 \quad \text{ for a constant}\,\, C.
\end{equation}

The function $\rho$ is unbounded on every interval $[t_0,t_*)$,
where $0\leq t_0<t_*$. Indeed, suppose that $0<\rho\leq L$ on such
an interval, and translate $t_0$ to $0$. Equation~\eqref{eq:vo-energy-bounds},
with $R_0=1+c_0L$, gives
\[
 |b|\rho\leq R_0,\qquad
 \rho\left(\sum_\alpha k_\alpha^2\right)^{1/2}\leq R_0.
\]
Cauchy--Schwarz therefore gives
\begin{equation}\label{eq:vo-weighted-d-bound}
 d\rho^2
 =\rho^2\sum_\alpha H_\alpha k_\alpha-(\rho f)^2
 \leq R_0\rho f-(\rho f)^2\leq\frac{R_0^2}{4}.
\end{equation}
By \eqref{eq:vo-Uderivative},
$d\rho^2\geq d(0)\rho(0)^2$. Together with \eqref{eq:vo-weighted-d-bound}, this gives
\[
 (\rho f)^2-R_0\rho f\leq-d(0)\rho(0)^2.
\]
Since $\rho f>0$, this quadratic inequality gives a uniform upper
bound for $\rho f$.

Suppose first that $d(t_1)>0$ for some $t_1\geq0$.
Equation~\eqref{eq:vo-Uderivative} gives
$d(t)\rho(t)^2\geq\varepsilon:=d(t_1)\rho(t_1)^2>0$ for $t\geq t_1$.
It follows from \eqref{eq:vo-weighted-d-bound} that
$\rho f\geq\varepsilon/R_0$. Substitution in
\eqref{eq:vo-Uderivative} yields
\[
 (d\rho^2)'=\frac{f'}{f}(\rho f)^2
       \geq\left(\frac{\varepsilon}{R_0}\right)^2\frac{f'}{f}.
\]
Integrating from $t_1$ to $t$ and using \eqref{eq:vo-weighted-d-bound},
we obtain
\[
 \left(\frac{\varepsilon}{R_0}\right)^2
 \log\frac{f(t)}{f(t_1)}
 \leq d(t)\rho(t)^2-d(t_1)\rho(t_1)^2
 \leq\frac{R_0^2}{4}-\varepsilon.
\]
Thus $f$ is bounded. The inequality $\rho f\geq\varepsilon/R_0$ then bounds
$\rho$ away from zero. Equations \eqref{eq:vo-first-integral} and
\eqref{eq:vo-second-integral} bound all the functions
$k_\alpha,b,H_\alpha,H_\alpha'$. If $t_*<\infty$, these bounds
allow continuation of the polynomial system beyond $t_*$, contrary
to maximality. If $t_*=\infty$, then
\eqref{eq:vo-positive} and $\rho\leq L$ imply $f'\geq\delta/L>0$,
contrary to boundedness of $f$.

It remains to consider $d\leq0$ on $[0,t_*)$. Define $s$ by
\begin{equation}\label{eq:vo-time-change}
 \frac{\mathrm ds}{\mathrm dt}=\frac1\rho,\qquad s(0)=0.
\end{equation}
A dot denotes differentiation with respect to $s$. By
\eqref{eq:vo-rho} and \eqref{eq:vo-time-change},
$\dot H_\alpha=\rho H_\alpha'$, $\dot f=\rho f'$, and
$\dot\rho=-b\rho^2$.
Since $d\leq0$ and $f^2-c>0$, \eqref{eq:vo-second-integral} gives
\begin{equation}\label{eq:vo-Hdot}
 \sum_\alpha\dot H_\alpha^2\leq C,
 \qquad \delta\leq\dot f\leq\sqrt C.
\end{equation}
In particular, $C\geq\delta^2>0$. Moreover,
\eqref{eq:vo-radial} gives
\begin{equation}\label{eq:vo-fddot}
 \ddot f=\rho^2(f''-bf')
       \geq2f(f^2-c)\rho^2\geq2f^3\rho^2.
\end{equation}

Let $s_*$ denote the right endpoint of the $s$-interval.
If $s_*<\infty$, integration of $|\dot\rho/\rho|=|b|\rho\leq R_0$ gives
\[
 \rho(s)\geq\rho(0)e^{-R_0s_*}>0,
 \qquad 0\leq s<s_*.
\]
Equations \eqref{eq:vo-energy-bounds} and \eqref{eq:vo-Hdot},
together with $f\leq f(0)+\sqrt{C}\,s_*$, bound
$k_\alpha,b,H_\alpha,H_\alpha'$. Since
\[
 t_* =\int_0^{s_*}\rho\,\mathrm ds\leq Ls_*<\infty,
\]
the original polynomial system extends beyond $t_*$, contrary to
maximality. Hence $s_*=\infty$.
Integrating \eqref{eq:vo-Hdot} and \eqref{eq:vo-fddot} yields
\begin{equation}\label{eq:vo-integrals}
 f(0)+\delta s\leq f(s)\leq f(0)+\sqrt{C}\,s,
 \qquad \int_0^\infty f^3\rho^2\,\mathrm ds<\infty.
\end{equation}

Since $f\geq f(0)>0$, equation~\eqref{eq:vo-integrals} gives
\[
 \int_0^\infty\rho^2\,\mathrm ds
 \leq\frac1{f(0)^3}\int_0^\infty f^3\rho^2\,\mathrm ds<\infty,
 \qquad
 \int_0^\infty\rho^2f^2\,\mathrm ds<\infty.
\]
Since $|H_\alpha|\leq f$, both $\rho$ and $\rho H_\alpha$ are square
integrable. By \eqref{eq:vo-rho} and \eqref{eq:vo-time-change},
\[
 \dot\rho=-(\rho b)\rho,
 \qquad
 \frac{\mathrm d}{\mathrm ds}(\rho H_\alpha)
 =\rho\dot H_\alpha-(\rho b)(\rho H_\alpha).
\]
The bounds for $\rho$, $\rho b$, $\rho f$, and $\dot H_\alpha$ therefore
give bounded derivatives of $\rho$ and $\rho H_\alpha$.
If one of these functions does not tend to zero, the derivative bound
gives disjoint intervals of a fixed positive length on which its absolute
value is bounded away from zero, contrary to square integrability.
It follows that
\begin{equation}\label{eq:vo-scaled-vanishing}
 \rho\to0,
 \qquad\rho H_\alpha\to0\quad\text{for every }\alpha.
\end{equation}

Differentiating $\rho b$ and $\rho k_\alpha$ and using
\eqref{eq:vo-G}, \eqref{eq:vo-b}, \eqref{eq:vo-rho}, and
\eqref{eq:vo-time-change}, we obtain
\begin{align}
 \frac{\mathrm d}{\mathrm ds}(\rho b)
 &=c\rho^2+\sum_\alpha
     \left(2(\rho H_\alpha)(\rho k_\alpha)-(\rho k_\alpha)^2\right),
 \label{eq:vo-b-product}\\
 \frac{\mathrm d}{\mathrm ds}(\rho k_\alpha)
 &=(\rho b)(\rho k_\alpha-2\rho H_\alpha).
 \label{eq:vo-k-product}
\end{align}
Since $c\leq0$, Young's inequality applied to \eqref{eq:vo-b-product} gives
\[
 \frac{\mathrm d}{\mathrm ds}(\rho b)
 \leq-\frac12\sum_\alpha(\rho k_\alpha)^2+2\rho^2f^2.
\]
The function
\[
 \rho b-2\int_0^s\rho(\sigma)^2f(\sigma)^2\,\mathrm d\sigma
\]
is nonincreasing and bounded below, because $|\rho b|\leq R_0$ and
\eqref{eq:vo-integrals} implies
$\int_0^\infty\rho^2f^2\,\mathrm ds<\infty$.
Consequently, $\rho b$ has a finite limit, and
$\int_0^\infty\sum_\alpha(\rho k_\alpha)^2\,\mathrm ds<\infty$.
By \eqref{eq:vo-energy-bounds} and \eqref{eq:vo-k-product}, the derivatives
of $\rho k_\alpha$ are bounded. Square integrability therefore yields
$\rho k_\alpha\to0$ for every $\alpha$. Substituting these limits and
\eqref{eq:vo-scaled-vanishing} into \eqref{eq:vo-first-integral}, we obtain
$E\geq0$ and $\rho b\to\pm\sqrt E$.
If this limit is negative, then $\rho b\leq-\varepsilon$ for some
$\varepsilon>0$ and all sufficiently large $s$. It follows that
$\dot\rho\geq\varepsilon\rho$, contrary to $\rho\to0$ in
\eqref{eq:vo-scaled-vanishing}.

If $E=0$, then \eqref{eq:vo-first-integral} gives
$|\rho b|\leq c_0\rho$. Thus
\[
 \frac{\mathrm d}{\mathrm ds}\left(\frac1\rho\right)=b\leq c_0,
 \qquad
 \rho(s)\geq\frac1{\rho(0)^{-1}+c_0s}.
\]
Together with $f(s)\geq f(0)+\delta s$, this contradicts the
integrability in \eqref{eq:vo-integrals}. Therefore
\begin{equation}\label{eq:vo-A-positive}
 E>0,\qquad \rho b\longrightarrow\sqrt E>0.
\end{equation}
Choose $a_0>0$ such that $\rho b\geq a_0$ for all sufficiently large $s$.
Equations \eqref{eq:vo-G}, \eqref{eq:vo-rho}, and
\eqref{eq:vo-time-change} give
\[
 \frac{\mathrm d}{\mathrm ds}
       \left(\rho^2(k_\alpha-H_\alpha)\right)
 =-\rho^2\dot H_\alpha.
\]
Since $\rho^2(k_\alpha-H_\alpha)\to0$ by
\eqref{eq:vo-scaled-vanishing} and $\rho k_\alpha\to0$,
integration from $s$ to infinity yields
\begin{equation}\label{eq:vo-backward}
 k_\alpha(s)-H_\alpha(s)
 =\int_s^\infty
       \left(\frac{\rho(\sigma)}{\rho(s)}\right)^2
       \dot H_\alpha(\sigma)\,\mathrm d\sigma.
\end{equation}
Since $\dot\rho/\rho=-\rho b$, equation~\eqref{eq:vo-A-positive} gives
$\rho(\sigma)/\rho(s)\leq e^{-a_0(\sigma-s)}$ for
$\sigma\geq s$ sufficiently large. Substituting this estimate and
\eqref{eq:vo-Hdot} into \eqref{eq:vo-backward}, we obtain
$|k_\alpha-H_\alpha|\leq\sqrt C/(2a_0)$.
Multiplying \eqref{eq:vo-tangent-product} by $\rho$ and using
\eqref{eq:vo-Hdot} and \eqref{eq:vo-integrals}, we conclude that
\[
 2\delta\bigl(f(0)+\delta s\bigr)
 \leq2f\dot f
 =\sum_\alpha(k_\alpha-H_\alpha)\dot H_\alpha
 \leq\frac{pC}{2a_0}
\]
for all sufficiently large $s$. This is impossible as $s\to\infty$.
Hence $\rho$ is unbounded on every interval $[t_0,t_*)$.

Equations \eqref{eq:vo-rho} and \eqref{eq:vo-energy-bounds} give
\begin{equation}\label{eq:vo-rho-growth}
 \rho'\leq1+c_0\rho,
 \qquad\rho(t)\leq e^{c_0t}\left(\rho(0)+t\right).
\end{equation}
If $t_*<\infty$, equation~\eqref{eq:vo-rho-growth} bounds $\rho$
on $[0,t_*)$, contrary to the unboundedness of $\rho$. Hence $t_*=\infty$.
Suppose that $f(t_0)>c_0$ for some $t_0\geq0$. By
\eqref{eq:vo-positive}, \eqref{eq:vo-Uderivative}, and
\eqref{eq:vo-energy-bounds}, we have
\[
 d(t_0)\rho(t_0)^2\leq d(t)\rho(t)^2
 \leq f(t)\left(\rho(t)-(f(t)-c_0)\rho(t)^2\right),
 \quad t\geq t_0.
\]
Since $f(t)\geq f(t_0)$, division by $f(t)$ yields
\[
 (f(t_0)-c_0)\rho(t)^2-\rho(t)
 \leq\frac{|d(t_0)|\rho(t_0)^2}{f(t_0)}.
\]
Since $f(t_0)-c_0>0$, the preceding inequality bounds $\rho$ on
$[t_0,\infty)$, a contradiction. Therefore
\begin{equation}\label{eq:vo-fbound}
 0<f(0)\leq f(t)\leq c_0,\qquad t\geq0.
\end{equation}
Equation~\eqref{eq:vo-fbound} is impossible if $c=0$.
For $c<0$, \eqref{eq:vo-q-positive} yields
\[
 (f'\rho)'\geq2c_0^2f(0)\rho.
\]
By \eqref{eq:vo-positive} and \eqref{eq:vo-fbound}, the function $f$
is increasing and bounded. Together with \eqref{eq:vo-rho-growth},
the preceding inequality contradicts Lemma~\ref{lem:growth-comparison},
applied with $\zeta=f$, $w=\rho$, and $b_0=2c_0^2f(0)>0$.
Therefore $H_\alpha\equiv0$ for every $\alpha$.
\end{proof}

\setlength{\bibsep}{6pt}


\begin{thebibliography}{99}

\bibitem{Abraham2026}
M. Abraham, \emph{Every biharmonic hypersurface in Euclidean space is minimal}, preprint, 2026.

\bibitem{AM2013}
K. Akutagawa and S. Maeta,
\emph{Biharmonic properly immersed submanifolds in Euclidean spaces},
Geom. Dedicata \textbf{164} (2013), 351--355.

\bibitem{AMOS2026}
\c{S}. Andronic, S. Montaldo, C. Oniciuc and A. Sanna,
\emph{Biconservative Weingarten surfaces with flat normal bundle in
$N^4(\epsilon)$}, Math. Nachr. \textbf{299} (2026), no.~7, 1624--1664.

\bibitem{BMO2008}
A. Balmu\c{s}, S. Montaldo and C. Oniciuc,
\emph{Classification results for biharmonic submanifolds in spheres},
Israel J. Math. \textbf{168} (2008), 201--220.

\bibitem{BMO2010}
A. Balmu\c{s}, S. Montaldo and C. Oniciuc,
\emph{Biharmonic hypersurfaces in 4-dimensional space forms},
Math. Nachr. \textbf{283} (2010), no.~12, 1696--1705.

\bibitem{BMO2013}
A. Balmu\c{s}, S. Montaldo and C. Oniciuc,
\emph{Biharmonic PNMC submanifolds in spheres},
Ark. Mat. \textbf{51} (2013), no.~2, 197--221.

\bibitem{BO2009}
A. Balmu\c{s} and C. Oniciuc,
\emph{Biharmonic surfaces of $\mathbb S^4$},
Kyushu J. Math. \textbf{63} (2009), no.~2, 339--345.

\bibitem{BO2012}
A. Balmu\c{s} and C. Oniciuc,
\emph{Biharmonic submanifolds with parallel mean curvature vector
field in spheres},
J. Math. Anal. Appl. \textbf{386} (2012), no.~2, 619--630.

\bibitem{BLO2021}
H. Bibi, E. Loubeau and C. Oniciuc,
\emph{Unique continuation property for biharmonic hypersurfaces
in spheres},
Ann. Global Anal. Geom. \textbf{60} (2021), no.~4, 807--827.

\bibitem{BrandingOniciuc2019}
V. Branding and C. Oniciuc,
\emph{Unique continuation theorems for biharmonic maps},
Bull. Lond. Math. Soc. \textbf{51} (2019), no.~4, 603--621.

\bibitem{CMO2001}
R. Caddeo, S. Montaldo and C. Oniciuc,
\emph{Biharmonic submanifolds of $\mathbb S^3$},
Internat. J. Math. \textbf{12} (2001), no.~8, 867--876.

\bibitem{CMO2002}
R. Caddeo, S. Montaldo and C. Oniciuc,
\emph{Biharmonic submanifolds in spheres},
Israel J. Math. \textbf{130} (2002), 109--123.

\bibitem{Chen1991}
B.-Y. Chen,
\emph{Some open problems and conjectures on submanifolds of finite type},
Soochow J. Math. \textbf{17} (1991), no.~2, 169--188.

\bibitem{Chen2019}
B.-Y. Chen,
\emph{Chen's biharmonic conjecture and submanifolds with parallel
normalized mean curvature vector},
Mathematics \textbf{7} (2019), no.~8, Paper No.~710, 7~pp.

\bibitem{Dimitric1992}
I. Dimitri\'c,
\emph{Submanifolds of $E^m$ with harmonic mean curvature vector},
Bull. Inst. Math. Acad. Sinica \textbf{20} (1992), no.~1, 53--65.

\bibitem{FLO2021}
D. Fetcu, E. Loubeau and C. Oniciuc,
\emph{Bochner--Simons formulas and the rigidity of biharmonic
submanifolds},
J. Geom. Anal. \textbf{31} (2021), no.~2, 1732--1755.

\bibitem{FetcuOniciuc2022}
D. Fetcu and C. Oniciuc,
\emph{Biharmonic and biconservative hypersurfaces in space forms},
Differential geometry and global analysis---in honor of Tadashi Nagano,
Contemp. Math., vol.~777, American Mathematical Society, Providence, RI,
2022, pp.~65--90.

\bibitem{FHT2026}
Y. Fu, M.-C. Hong and G. Tian,
\emph{The biharmonic hypersurface flow and the Willmore flow in
higher dimensions},
J. Eur. Math. Soc. (2026), online first,
\href{https://doi.org/10.4171/JEMS/1770}{doi:10.4171/JEMS/1770}.

\bibitem{FHY2026}
Y. Fu, M.-C. Hong and D. Yang,
\emph{Biharmonic hypersurfaces in space forms},
arXiv:2610.04165v2 (2026).

\bibitem{FHZ2021}
Y. Fu, M.-C. Hong and X. Zhan,
\emph{On Chen's biharmonic conjecture for hypersurfaces in $\mathbb R^5$},
Adv. Math. \textbf{383} (2021), Paper No.~107697, 28~pp.

\bibitem{FHZ2023}
Y. Fu, M.-C. Hong and X. Zhan,
\emph{Biharmonic conjectures on hypersurfaces in a space form},
Trans. Amer. Math. Soc. \textbf{376} (2023), no.~12, 8411--8445.

\bibitem{HV1995}
Th. Hasanis and Th. Vlachos,
\emph{Hypersurfaces in $E^4$ with harmonic mean curvature vector field},
Math. Nachr. \textbf{172} (1995), no.~1, 145--169.

\bibitem{Jiang1986}
G. Y. Jiang,
\emph{2-harmonic maps and their first and second variational formulas}
(Chinese), Chinese Ann. Math. Ser. A \textbf{7} (1986), no.~4, 389--402.
English summary, Chinese Ann. Math. Ser. B \textbf{7} (1986), no.~4, 523.

\bibitem{KN1969}
S. Kobayashi and K. Nomizu,
\emph{Foundations of Differential Geometry}, Vol.~II,
Interscience Publishers, New York--London--Sydney, 1969.

\bibitem{Maeta2014}
S. Maeta,
\emph{Properly immersed submanifolds in complete Riemannian manifolds},
Adv. Math. \textbf{253} (2014), 139--151.

\bibitem{MaetaWintgen2026}
S. Maeta,
\emph{Biharmonic Wintgen ideal submanifolds in Riemannian manifolds
of constant sectional curvature}, arXiv:2606.25791v2 (2026).

\bibitem{MaetaPseudo2026}
S. Maeta,
\emph{Four-di\-men\-sion\-al pseudo-um\-bil\-i\-cal biharmonic submanifolds
in six-di\-men\-sion\-al Riemannian manifolds of constant sectional
curvature}, arXiv:\allowbreak 2609.17136v2 (2026).

\bibitem{NU2013}
N. Nakauchi and H. Urakawa,
\emph{Biharmonic submanifolds in a Riemannian manifold with non-positive
curvature}, Results Math. \textbf{63} (2013), no.~1--2, 467--474.

\bibitem{OuChen2020}
Y.-L. Ou and B.-Y. Chen,
\emph{Biharmonic submanifolds and biharmonic maps in Riemannian geometry},
World Scientific Publishing, Hackensack, NJ, 2020.

\bibitem{OuTang2012}
Y.-L. Ou and L. Tang,
\emph{On the generalized Chen's conjecture on biharmonic submanifolds},
Michigan Math. J. \textbf{61} (2012), no.~3, 531--542.

\bibitem{ST2023}
R. \c{S}en and N. C. Turgay,
\emph{Biharmonic PNMCV submanifolds in Euclidean 5-space},
Turkish J. Math. \textbf{47} (2023), no.~1, 296--316.

\end{thebibliography}
\end{document}